\documentclass[a4paper, 12pt]{article}
\usepackage[utf8]{inputenc}
\usepackage{amssymb}\usepackage{multirow}
\usepackage{bbold}
\usepackage{xcolor}
\usepackage[T1]{fontenc}
\usepackage[utf8]{luainputenc}
\usepackage{amsmath}
\usepackage{bbm}
\usepackage{import}
\usepackage[english]{babel}
\usepackage{blindtext}
\usepackage[export]{adjustbox}
\usepackage{subcaption}
\usepackage{amssymb}
\usepackage{amsthm}
\usepackage{float}
\usepackage[ruled,vlined]{algorithm2e}
\usepackage[margin=2cm]{geometry}
\usepackage[authoryear]{natbib}
\usepackage{xcolor}
\usepackage{amsthm}
\usepackage{caption}
\usepackage{scalerel,stackengine}
\stackMath
\newcommand{\cov}{\text{Cov}}
\newcommand{\var}{\text{Var}}
\newcommand\reallywidehat[1]{%
\savestack{\tmpbox}{\stretchto{%
  \scaleto{%
    \scalerel*[\widthof{\ensuremath{#1}}]{\kern-.6pt\bigwedge\kern-.6pt}%
    {\rule[-\textheight/2]{1ex}{\textheight}}
  }{\textheight}%
}{0.5ex}}%
\stackon[1pt]{#1}{\tmpbox}%
}
\newtheorem{definition}{Definition}
\newtheorem{proposition}{Proposition}[]

\usepackage{url}
\usepackage{enumitem}
\setlist[enumerate]{itemsep=0mm}
\usepackage{outlines}
\usepackage{yhmath}
\usepackage{amsthm,amsmath,amsfonts,amssymb}
\RequirePackage{mathtools}

\usepackage{xcolor}
\usepackage{bm,accents}
\usepackage{graphicx}
\usepackage{nicefrac}
\usepackage[super]{nth}

\DeclareMathOperator*{\argmin}{arg\,min}
\usepackage{color}
\bibpunct{(}{)}{;}{a}{,}{,}
\def\spacingset#1{\renewcommand{\baselinestretch}%
	{#1}\small\normalsize} \spacingset{1.6}
\renewcommand{\baselinestretch}{1.4} 

\newtheorem{corollary}{Corollary}[]

\theoremstyle{remark}
\newtheorem{remark}{Remark}[]
\usepackage{authblk}

\title{A zero-estimator approach for estimating the signal level in a    high-dimensional model-free  setting}
\author{Ilan Livne, David Azriel,  Yair Goldberg}
\affil{Technion - Israel Institute of Technology}
\date{\today}

\begin{document}
\maketitle

\begin{abstract}
We study a high-dimensional regression setting under the assumption of known covariate distribution.  We aim at estimating the amount of  explained variation in the response  by the best linear function of the covariates  (the signal level). In our setting, neither sparsity of the coefficient vector, nor normality of the covariates or linearity of the conditional expectation are assumed.
We present  an unbiased and consistent estimator and then 
 improve it by using a zero-estimator approach, where a zero-estimator is a statistic whose expected value is zero.
 More generally, we present an algorithm based on the zero estimator approach that in principle can improve any given estimator. 
We study some asymptotic properties of the proposed estimators and demonstrate their finite sample performance in a simulation study.

\vspace{9pt}
\noindent {\it Key words and phrases:}
Linear Projection, Semi-supervised setting, U-statistics,
Variance estimation, Zero-estimators.
  
\end{abstract}

\section{Introduction}
In many regression settings, an important goal is to estimate the signal and noise levels, i.e., to quantify the amount of variance in the response variable that  can be explained  by a set of covariates, versus how much of the variation is left unexplained.  
When the covariates' dimension is low and a linear regression model is assumed, the ordinary least squares method can be used to find  a consistent estimator for the signal level.
However, in a high-dimensional setting, the  least squares method breaks down   and it becomes more challenging to develop good estimators without further assumptions.
In recent years, several methods have been proposed for estimating the signal level  under the assumption that the regression coefficient vector $\beta$ is sparse 
(\citet{fan2012variance,sun2012scaled,chatterjee2015prediction,verzelen2018adaptive,tony2020semisupervised}).
Other widely-used methods
assume some probabilistic structure on~$\beta$  (e.g., $\beta$ is Gaussian) and  use  maximum likelihood  to derive consistent estimators of the signal level 
(\citet{yang2010common,bonnet2015heritability}).
  These methods have been extensively studied  in the literature of \emph{random-effect models} where $\beta$ is treated as random; see \citet{de2015genomic} and references therein.
However,  methods that rely on the assumption that   $\beta$ is either  sparse or highly structured  may not perform well when these assumptions fail to hold.  
For example, a known problem in genetics is the problem of missing heritability (\citet{de2015genomic,zhu2020statistical}). 
   Heritability is  defined as  the fraction of the observed outcome (phenotype) that is explained by genetic factors. 
    The term ``missing heritability'' is traditionally used to describe the gap between heritability estimates from genome-wide-association-studies (GWAS) and the corresponding estimates from family studies. To explain the gap, it has been suggested that some phenotypes are explained by a numerous number of genetic factors that their individual effect is too small to detect, but their collective effect is significant (\citet{yang2010common, young2022discovering}). In such a setting,   methods that rely on the sparsity assumption may fail to provide accurate estimates.

Rather than assuming sparsity or other structural assumptions on $\beta$, a different approach for estimating the signal level in a high-dimensional setting  is  to assume some or complete knowledge  about the covariate distribution. This can be justified, for example, in the semi-supervised setting when one has access to a large amount of unlabeled (covariate) data  without the corresponding labels (responses).  
 When the covariates are assumed independent Gaussian, \citet{Dicker}  proposed estimators
 based on the method-of-moments  and \citet{janson2017eigenprism}  used convex optimization techniques.
    In both methods, the Gaussian assumption was used to  show  consistency and asymptotic-normality,   and it is not clear how  robust  these methods are  when  the assumptions are violated.
    Dropping the Gaussian independent covariate assumption,
\citet{livne2021improved}
 proposed a consistent estimator  
 under the assumption  that the first two moments of covariates are known. 
 More recently, \citet{chen2022statistical}
 proposed an estimator that 
   is consistent and asymptotically-normal when the covariates are independent and  the entries of $\beta$ are small and dense.

All of the estimators that we reviewed above were developed under the assumption that the linear model is true, which can be unrealistic in many situations.
In this work, we focus our attention on the \emph{model-free} setting, i.e., no assumptions are made about the relationship between the covariates and the response.
Under this setting,  \citet{kong2018estimating}
proposed a   consistent
 estimator under some assumptions on the covariance matrix. 
In this paper we follow the two-stage approach  presented in \citet{livne2021improved},  where an initial estimator is first suggested and then a  zero-estimator  is used to reduce its variance. 
Our initial estimator is the same as in  \citet{kong2018estimating} and 
\citet{livne2021improved}, and the zero estimators we use are tailored to the model-free framework. Furthermore, we provide a general algorithm that, in principle, improve any initial estimator and we also demonstrate the usefulness of the algorithm for several initial estimators.

 The rest of this work is organized as follows. 
In Section \ref{sec:Preliminaries}, 
 we discuss the parameters of interest in a  model-free setting under the assumption that the first two  moments of the covariates are known.
  In Section \ref{sec: naive_estimators}, we present our initial estimators and prove that they are consistent under some minimal assumptions.
 In Section \ref{sec:improv_of_Naive_sec}, we use the zero-estimator approach  to construct two improved estimators and then study  some theoretical properties of the improved estimators. Simulation results are  given in Section \ref{sim_res}. Section \ref{gener_es} demonstrates  how the zero-estimator approach can be  generalized to other estimators. A discussion is given in Section \ref{discuss}. The proofs are provided in the Appendix.

\section{\label{sec:Preliminaries}Preliminaries}
Let $X\in {{\mathbb{R}}^{p}}$ be a random vector of covariates and let $Y \in \mathbb{R}$ be the response. 
The conditional mean $E(Y|X)$ is the \emph{best predictor} in the sense that it minimizes 
 the mean squared error $E\{[Y-g(X)]^2\}$ over all measurable functions $g(X)$  (see \citet{hansen_2022}, p.~$25$).
However, the functional form of $E(Y|X)$ is typically unknown and difficult to estimate, especially in a high-dimensional setting. Consequently, we can define the best linear approximation to $E(Y|X)$. 
\begin{definition}\label{def:BLP}
Assume that both $E(Y^2)$, \(E( {{{\left\| X \right\|}^2}} )\) exist and that  the covariance matrix of $X$, denoted by  ${\Sigma_X}$,
  is invertible. Then, the \emph{best linear predictor}, $\alpha+\beta^TX$, is defined by the unique $\alpha$ and $\beta$ that minimize the mean squared  error 
 	$$(\alpha,\beta)=\argmin_{a\in\mathbb{R},b\in\mathbb{R}^p}E(Y-a-b^TX)^2,$$
and, by \citet{hansen_2022}, pp.~34-36, is given by   \begin{equation}\label{BLP_def}
 \beta={\Sigma_X}^{-1}\{E(XY)-E(X)E(Y)\}\quad\text{and}\quad  \alpha=E(Y)-\beta^TE(X).   
\end{equation}
\end{definition}
The best linear predictor is essentially the population version of the OLS method. Notice that $\alpha,\beta$ also satisfy  $(\alpha,\beta)=\argmin_{a\in\mathbb{R},b\in\mathbb{R}^p}E[E(Y|X)-a-b^TX]^2$.
It is a model-free quantity, i.e.,  no specific assumptions are made about the relationship between $X$ and $Y$. In particular, we do not  assume that $E(Y|X)$ is linear in $X$. 
If $E(Y|X)$ happens to be linear, say $E(Y|X)=\tilde{\alpha}+\tilde{\beta}^T X$, then  the best linear predictor parameters $(\alpha,\beta)$ coincide with the model parameters $(\tilde{\alpha},\tilde{\beta}).$
However, when $E(Y|X)$ is not linear, the parameters $\alpha$ and $\beta$ are still meaningful: they describe the overall direction of the association between $X$ and $Y$ (\citet{buja2019models}).
Hence, Definition \ref{def:BLP} is  useful
since in most cases we have no reason to believe  that $E(Y|X)$ is indeed linear in $X$.

We now wish to decompose the variance of $Y$ into signal and noise levels. Let $\sigma_Y^2$ denote the variance of $Y$ and define the residual $\epsilon\equiv Y-\alpha-\beta^TX$. Notice that both  $E(X\epsilon)= 0$ and $E(\epsilon)=0$ by construction.  Write
\begin{equation}\label{eq:var y decompose}
	\sigma_{Y}^2 =\text{Var}(\alpha+X^T\beta+\epsilon)
	 =\beta^T\text{Var}(X)\beta+\text{Var}(\epsilon)
	 =\beta^T{\Sigma}_X\beta+\sigma^2,
	 \end{equation}
where $\sigma^2\equiv\var(\epsilon)$ and ${\Sigma}_X\equiv\var(X)$. 
Here, the \emph{signal} level $\tau^2\equiv\beta^T{\Sigma}_X\beta$ can be thought of as the total variance explained by the best linear function of the covariates. The \emph{noise} level~$\sigma^2$ is the variance left unexplained.
Notice that the parameters $\tau^2$, $\sigma^2$ and $\sigma_Y^2$ depend on $p$ but this is suppressed in the notation.

A common starting point for many works of regression problems is to use strong assumptions about $E(Y|X)$, and minimal assumptions, if any,  about the covariate $X$. In this work, we take the opposite approach: we make no assumptions about $E(Y|X)$ but assume we know everything about $X$ instead.
This can be justified, for example, in a  semi-supervised setting where, in addition to the labeled data, we also have access to a large amount of unlabeled data; see, for example, the work of \citet{zhang2019high} who study estimation of the mean and variance of $Y$ in a semi-supervised, model-free setting.
 Note that the setting of known-covariate distribution has already been presented and discussed in the context of high-dimensional regression   as in \citet{candes2018panning,berrett2020conditional} and \citet{wang2020power}. 
 Our  goal here is to develop good estimators for $\tau^2$ and $\sigma^2$ in a model-free setting under the assumption that the distribution of the covariates is known.
 
 Let $(X_1,Y_1),...,(X_n,Y_n)$ be i.i.d.\ observations drawn from an unknown distribution  where $X_i\in\mathbb{R}^p$ and $Y_i\in\mathbb{R} $.  Let $(X,Y)$ denote a generic observation from the sample.
We assume that $E(X)\equiv \mu$ is known and also that the variance matrix ${\Sigma}_X$ is known and invertible. Linear transformations do not affect the signal and noise levels. Thus, we can apply the transformation
$X\mapsto{\Sigma}_X^{-1/2}(X-~\mu)$ and assume w.l.o.g. that
\begin{equation}\label{assum}
E(X)=0\quad \text{and} \quad  {\Sigma}_X=\textbf{I}. \end{equation}
By  \eqref{eq:var y decompose},  
 $\sigma_{Y}^2=\|\beta\|^{2}+\sigma^2$, 
  which means that in order to estimate $\sigma^2$ it is enough to estimate both $\sigma_{Y}^2$ and $\tau^2\equiv\|\beta\|^{2}$. The variance term  $\sigma_{Y}^2$ can be easily estimated from the sample. Hence,  the main challenge is to derive an estimator for $\tau^2.$

\section{Initial Estimators}\label{sec: naive_estimators}
In this section, we present our initial estimators for the signal and noise levels,  $\tau^2$ and $\sigma^2$.
Interestingly, when the linear model is true and $X$ is assumed to be constant, no consistent estimator of $\sigma^2$ exists (\citet{azriel2019conditionality}).
However,  when 
 $X$ is random, a consistent estimator does exist if ${\Sigma}$ is known (\citet{verzelen2018adaptive}).
 The current  work goes one step further as we generalize this result without assuming linearity.
 Indeed, Proposition \ref{consistency_naive}  below 
 demonstrates that by knowing the first and second moments of $X$, it is possible, under some mild assumptions, to construct  consistent estimators of $\tau^2$ and $\sigma^2$ in a high-dimensional setting without assuming that  $E(Y|X)$ is linear. The estimator we use below was suggested by \citet{kong2018estimating} who provided an upper bound on the variance. Our analysis is more general and we discuss sufficient conditions for consistency.
 
 Let $W_{ij}\equiv X_{ij}Y_i$ for $i=1,...,n$, and $j=1,...,p$. Notice that
\begin{equation}\label{beta}
  E\left( {{W_{ij}}} \right) = E\left( {{X_{ij}}{Y_i}} \right) = E\left[ {{X_{ij}}\left( {{\alpha + \beta ^T}{X_i} + {\varepsilon _i}} \right)} \right] = {\beta _j},  
\end{equation}
where  in the last equality we used \eqref{assum} and the orthogonality between $X$ and $\epsilon.$ 
   Now, since
   $\{E(W_{ij})\}^2=E(W_{ij}^2)-\text{Var}(W_{ij})$, a natural unbiased estimator for $\beta_j^2$ is
 \begin{equation}\label{beta_j_hat} 
 {\hat\beta_j^2}\equiv\frac{1}{n}\sum\limits_{i=1}^{n}W_{ij}^2-\frac{1}{n-1}\sum\limits_{i=1}^{n}(W_{ij}-\overline{W}_j)^2 = \binom{n}{2}^{-1}\sum_{i_1< i_2}^n W_{i_1j}W_{i_2j}, 
 \end{equation}
 where $\overline{W}_j=\frac{1}{n}\sum_{i=1}^{n}W_{ij}$. Thus, an unbiased estimator of $\tau^2$  is given by
\begin{equation} \label{estimates}
{\hat \tau ^2} = \sum\limits_{j = 1}^p {\hat \beta _j^2 = \binom{n}{2}^{-1}\sum\limits_{{i_1} < {i_2}}^{} {W_{{i_1}}^T{W_{{i_2}}}} },
\end{equation}
where \({W_i} = \left( {{W_{i1}},...,{W_{ip}}} \right)^T\).
We call  $\hat\tau^2$  the \emph{naive estimator}.
Notice that $\hat\tau^2$ is a U-statistic with the kernel 
\(h\left( {{W_1},{W_2}} \right) = W_1^T{W_2}\) and thus its variance can be calculated directly by using U-statistic properties  (see \citet{van2000asymptotic}, Theorem 12.3).

Let 
\({\zeta _1} \equiv {\beta ^T}{\bf{A}}\beta  - {\left\| \beta  \right\|^4}\), \({\zeta _2} \equiv \left\| {\bf{A}} \right\|_F^2 - {\left\| \beta  \right\|^4}\),
\({\bf{A}} \equiv E\left( {{W}W^T} \right)\)
and $\|\mathbf{A}\|_F$ denotes the Frobenius norm of $\bf{A}.$ The following proposition calculates the variance of $\hat{\tau}^2$.
\begin{proposition}\label{prop:var_naive_tau2}
Assuming that $\zeta_1$ and $\zeta_2$ are finite,
the variance of the naive estimator $\hat\tau^2$ is given by
\begin{equation}\label{eq:var_naive}
\var\left( {{{\hat \tau }^2}} \right) = \frac{{4\left( {n - 2} \right)}}{{n\left( {n - 1} \right)}}{\zeta _1} + \frac{2}{{n\left( {n - 1} \right)}}{\zeta _2}.
\end{equation}
\end{proposition}
The next proposition shows that the naive estimator $\hat\tau^2$ is consistent under some  assumptions.  
\begin{proposition}\label{consistency_naive}
Assume that $\frac{\|\textbf{A}\|^{2}_F}{n^2}\xrightarrow{n\rightarrow\infty}0,$ and that $\tau^2$ is bounded. Then,
$$\hat\tau^2 - \tau^2 \overset{p}{\rightarrow}~0.$$
 \end{proposition}
 Similarly, an  estimator for the noise level $\sigma^2$ can be obtained by 
\begin{equation}\label{eq: sigma2_estimator}
{\hat \sigma ^2} \equiv \hat \sigma _Y^2 - {\hat \tau ^2},    
\end{equation}
 where 
\(\hat \sigma _Y^2 \equiv \frac{1}{{n - 1}}\sum\limits_{i = 1}^n {{{\left( {{Y_i} - \bar Y} \right)}^2}} \) is the standard unbiased estimator of $\sigma_Y^2.$
Let \({\mu _4} \equiv E\left[ {{{\left( {Y - \alpha } \right)}^4}} \right]\), \(\pi  = {\left( {{\pi _1},...,{\pi _p}} \right)^T}\) where \({\pi _j} \equiv E\left[ {{{\left( {Y_1 - \alpha } \right)}^2}{W_{1j}}} \right]\).
The variance of $\hat\sigma^2$ is given by the following proposition.
\begin{proposition}\label{proposition: consistent_sigma2}
The variance of $\hat\sigma^2$ is
 \begin{align}\label{eq:var_sigma_2_hat}
\var\left( {{{\hat \sigma }^2}} \right) &= \left[ {\frac{1}{n}{\mu _4} - \frac{{\left( {n - 3} \right)}}{{n\left( {n - 1} \right)}}\sigma _Y^4} \right] + \var(\hat\tau^2)\nonumber\\
&- \frac{4}{n}\left( {{\pi ^T}\beta  - {\tau ^2}\sigma _Y^2} \right) + \frac{4}{{n\left( {n - 1} \right)}}\sum\limits_{j = 1}^p {{{\left\{ {E\left[ {{W_{1j}}\left( {Y_1 - \alpha } \right)} \right]} \right\}}^2}},
\end{align}
where $\var(\hat\tau^2)$ is given in \eqref{eq:var_naive}.
\end{proposition}
The following result is a corollary of Propositions 2 and 3. 
\begin{corollary}\label{cor:consistency_naive_sigma2}
Assume that $\mu_4$ and $\tau^2$  are bounded and that $\frac{\|{\bf A}\|_F^{2}}{n^2}\rightarrow 0.$ Then,
$$\hat\sigma^2-\sigma^2\overset{p}{\rightarrow}~0.$$
\end{corollary}
 The condition $\frac{\|{\bf A}\|_F^{2}}{n^2}\rightarrow 0$   holds in various settings. For example, it can be shown to hold  when $p/n^2\rightarrow 0 $  and $Y$ is bounded.  For more   examples and details, see Remark \ref{remark: forb_A} in the Appendix.

\section{Reducing Variance Using a Zero Estimator }\label{sec:improv_of_Naive_sec}
In this section, we study how  the naive estimator $\hat\tau^2$, and consequently  $\hat\sigma^2$, can be improved by
 using the
assumption that the distribution of $X$ is known.
  We use  zero-estimators to construct an improved  unbiased estimator of $\tau^2$.
This is also   known as the method of \emph{control variables} from the Monte-Carlo literature; see, e.g., \citet{glynn2002some,lavenberg1981perspective}.
  Here, a zero-estimator is defined as a statistic $Z$ such that $E(Z)=0$.  For a given zero-estimator $Z$ and a constant $c$, we define a new estimator $\hat\tau^2(Z,c)$ as
  \begin{equation}\label{improves_est}
   \hat\tau^2(Z,c)=\hat\tau^2-cZ.
  \end{equation}
  For a fixed $c$, notice that  $\hat\tau^2(Z,c)$ is an unbiased estimator  for $\tau^2.$
   Also notice that for every function $f$ of the covariates $X_1,...,X_n$,  one can always define a zero-estimator  $${Z} =  {f\left( {{X_1,\dots,X_n}} \right) - E\left[ {f\left( {{X_1,\dots,X_n}} \right)} \right]}.$$ This is possible since we assume that the distribution of the covariates is known and hence $E[f\left( {{X_1,\dots,X_n}} \right)]$ is known.
The variance of $\hat\tau^2(Z,c)$ is
$$\var[\hat\tau^2(Z,c)]= \var\left( {{{\hat \tau }^2}} \right) + {c^2}\var\left( {{ Z}} \right) - 2c\,\cov\left( {{{\hat \tau }^2},{ Z}} \right).$$
Minimizing the variance with respect to  $c$ yields the minimizer 
\begin{equation}\label{eq:c_star}
{c^*} = \frac{{\cov\left( {{{\hat \tau }^2},{ Z}} \right)}}{{\var\left( {{ Z}} \right)}}.    
\end{equation}
Hence, the corresponding oracle-estimator  is  $$\hat\tau^2(Z,c^*)=\hat\tau^2-c^*Z.$$ 
We use the term \emph{oracle} since the optimal coefficient $c^*$ is an unknown quantity.
The  variance of the above oracle-estimator is 
\begin{align}\label{MSE_T_c_star}
 \var[\hat\tau^2(Z,c^*)] &= \var\left( {{{\hat \tau }^2}} \right) - \frac{{\left[ {\cov\left( {{{\hat \tau }^2},{ Z}} \right)} \right]^2}}{{\var\left( {{ Z}} \right)}} \nonumber\\
&= \left( {1 - \rho _{{{\hat \tau }^2},{Z}}^2} \right)\var\left( {{{\hat \tau }^2}} \right),   
\end{align}
where $\rho _{{{\hat \tau }^2},{Z}}$ is the correlation coefficient between $\hat\tau^2$ and $Z$.
The term ${1 - \rho _{{{\hat \tau }^2},{Z}}^2}$ is the factor by which  $\var(\hat\tau^2)$ could be reduced if the optimal  coefficient $c^*$ was known.
Thus, the more correlation there is between the zero-estimator $Z$ and the naive estimator $\hat\tau^2$, the greater the reduction in variance.

There are two challenges  to be addressed with the above approach. First, one should  find a simple zero-estimator $Z$ which is  correlated with the naive estimator $\hat\tau^2$. Second, the optimal  coefficient~$c^*$ is an unknown quantity and therefore  needs to be estimated.

To address the first challenge,  we propose the following  zero-estimator 
 $$Z_g\equiv{\bar g_{n}} = \frac{1}{n}\sum\limits_{i = 1}^n {{g}(X_i)},$$ where
    ${g}(X_i) = \sum\limits_{j < j'}^{} {{X_{ij}}{X_{ij'}}}.$ 
    In Remark \ref{c_star_single} in the  Appendix, we show that the optimal coefficient, with respect to $Z_g$, is
\begin{equation}\label{eq:c_star_g}
{ c_g^*} = \frac{{{2\beta ^T}\theta_g }}{{\var [ {{g(X)}}] }},
\end{equation} 
where \({\theta _g} = E\left[ {Wg\left( X \right)} \right].\)
Notice that $\var [{g(X)}]$   is a  known quantity since the distribution of $g(X)$  is assumed to be known. 
Hence, the corresponding oracle-estimator  is 
\begin{align}\label{eq:T_g}
    T_g\equiv \hat\tau^2(Z_g,c_g^*)=\hat\tau^2-c_g^*Z_g.
\end{align}

To address the second challenge, i.e., to estimate the optimal coefficient $c_g^*$, we suggest the following  unbiased U-statistic  estimator
\begin{equation}\label{eq:c_g_hat_star}
\hat c_g^* =  \frac{\binom{n}{2}^{-1}{\sum\limits_{{i_1} \neq {i_2}}^{} {W_{{i_1}}^T{W_{{i_2}}}g\left( {{X_{{i_2}}}} \right)} }}{{\var\left[ {g\left( X \right)} \right]}}.
\end{equation}
Thus,  the corresponding  improved  estimator is \begin{equation}\label{eq:T_g_hat}
T_{\hat g}\equiv\hat\tau^2(Z_g,\hat c_g^*)= {\hat \tau ^2} - {\hat c_g^*}Z_g.
\end{equation}
Using the zero-estimator $Z_g$ has a potential drawback. It uses \emph{all}   the $p$ covariates of the vector~$X$   regardless of the sparsity level in the data, which can result in some additional variability due to unnecessary estimation.  Intuitively, when the sparsity level is high, i.e., only a small number of covariates  plays an important role in explaining the response $Y$,  
it is  inefficient to use a zero-estimator that incorporates all the $p$ covariates.
In such a setting, it is reasonable to modify the zero-estimator $Z_g$ such that only
 a small set of covariates  will be included, preferably the  covariates  that capture  a significant part of the signal level $\tau^2.$
  Selecting such a set of covariates can be difficult and one may use a covariate-selection procedure for this purpose.

We call $\delta$ a covariate \emph{selection procedure} if for every dataset   it chooses a subset of indices \({{\bf{S}}_\delta }\) from \(\left\{ {1,...,p} \right\}\).
 Different covariate-selection methods exist in the literature (see  \citet{oda2020fast} and references therein) but  these are not a primary focus of this work.
  For a given selection procedure $\delta$ we  modify the estimator $T_{\hat g}$ such that only  the indices in~${\bf{S}}_\delta$ will be included in its zero-estimator term.
   This modified estimator, which is based on a given selection procedure~$\delta$, is presented in the algorithm below. 
  
\vspace{10mm}
\begin{algorithm}[H]\label{alg:1}
 \caption{Proposed estimator for  $\tau^2.$ }
\SetAlgoLined

\vspace{0.4 cm}

\textbf{Input:}
 A dataset \(\left( {{{\bf{X}}_{n \times p}},{Y_{n \times 1}}} \right)\) and a selection procedure $\delta$.
 \begin{enumerate}
 \item Calculate the naive estimator 
 ${\hat \tau ^2} =  \binom{n}{2}^{-1}\sum\limits_{{i_1} < {i_2}}^{} {W_{{i_1}}^T{W_{{i_2}}} },$
          \item Apply procedure $\delta$ to  \(\left( {{{\bf{X}}_{n \times p}},{Y_{n \times 1}}} \right)\) to construct the set \({{\bf{S}}_{\delta}}.\) 
      \end{enumerate}
\textbf{Output}:
Return the estimator
\begin{equation}\label{selec_est}
 T_{\hat h}\equiv  {\hat \tau ^2} - \hat c_h^*{Z_h},   
\end{equation}
where  
\(\hat c_h^* = \frac{{\binom{n}{2}^{-1} \sum\limits_{{i_1} \neq {i_2}}^{} {W_{{i_1}}^T{W_{{i_2}}}h\left( {{X_{{i_2}}}} \right)} }}{{\var \left[ {h\left( {{X}} \right)} \right]}}\), 
 ${Z_h} = \frac{1}{n}\sum\limits_{i = 1}^n {{h}\left( {{X_i}} \right)} $ and  \(h\left( {{X_i}} \right) = \sum\limits_{j < j' \in {\bf{S}}_\delta}^{} {{X_{ij}}{X_{ij'}}} \).
\end{algorithm}
\vspace{10mm}

Notice that the estimator $T_{\hat g}$ 
defined in \eqref{eq:T_g_hat} is a special case of the estimator $T_{\hat h}$
defined in Algorithm~\ref{alg:1},   when ${\bf{S}}_\delta = \{1,\dots p \}$, i.e., when  $\delta$ selects all the $p$ covariates.

Recall  that in this work we  treat $p$ as a function of $n$, i.e., $p\equiv p_n$ but this is suppressed
in the notation.
Let  $\Theta \equiv \Theta_n \subseteq \{1, \ldots,p\}$  be a deterministic sequence of subsets.
In order to  analyze the estimator 
$T_{\hat h}$
 we define a stability  property, which is given next.
 \begin{definition}\label{def:stable}
A selection procedure $\delta$ is \emph{stable}  
if there exists a deterministic sequence of subsets $\Theta$  such that
\begin{equation}\label{stable_condition}
\lim_{n\to \infty} n[P( {\bf S}_\delta  \ne  \Theta)]^{1/2}  = 0. 
\end{equation}
\end{definition}
Definition \ref{def:stable} states  that  a selection procedure $\delta$ is stable if it is asymptotically close to a deterministic procedure at a suitable rate.
The convergence rate of  many practical selection procedures
is exponential,
which is much faster than is required for the condition  to hold. For example,  the lasso algorithm  asymptotically selects the support of $\beta$ at an exponential rate   under some assumptions
 (see \citet{hastie_tibshirani_wainwright_2015}, Theorem 11.3). Notice also that the stability condition holds trivially when ${\bf{S}}_\delta = \{1,\dots p \}$, i.e., when $\delta$  selects all the $p$ covariates for all $n$.
 
 Define the oracle-estimator ${T_{h}} \equiv {\hat \tau ^2}\left( {{Z_h}, c_h^*} \right)$, where \(c_h^* = \frac{{{\beta ^T}{\theta _h}}}{{\var\left[ {h\left( X \right)} \right]}}\) and \({\theta _h} = E\left[ {Wh\left( X \right)} \right]\).
  Let
\begin{align}\label{eq:list_def}
f\left( {{X_i}} \right) &\equiv \sum\limits_{j < j' \in \Theta}^{} {{X_{ij}}{X_{ij'}}}, 
\quad
{T_{ f}} \equiv {\hat \tau ^2}\left( {{Z_f}, 
c_f^*} \right),
\quad
{T_{\hat f}} \equiv {\hat \tau ^2}\left( {{Z_f},\hat c_f^*} \right), 
\quad
c_f^* \equiv \frac{{2{\beta ^T}{\theta _f}}}{{\var\left[ {f\left( X \right)} \right]}},\nonumber\\
 \hat c_f^* &\equiv \frac{{\frac{2}{{n\left( {n - 1} \right)}}\sum\limits_{{i_1} \ne {i_2}}^{} {W_{{i_1}}^T{W_{{i_2}}}f\left( {{X_{{i_2}}}} \right)} }}{{\var\left[ {f\left( X \right)} \right]}},
 \quad
  \theta_f\equiv E[Wf(X)],
  \quad
  {Z_f} \equiv \frac{1}{n}\sum\limits_{i = 1}^n {f\left( {{X_i}} \right)},\nonumber\\
    {\bf{B}} &= E\left[ {W{W^T}f\left( X \right)} \right], \quad   {\bf{C}} = E\left[ {W{W^T}f^2{{\left( X \right)}}} \right].  
  \end{align} 
  
   We now prove that the proposed estimator in Algorithm \ref{alg:1} is asymptotically equivalent to its oracle version $T_h$ under some conditions.
  \begin{proposition}\label{prop:singel_asymptotic}
Assume that the selection procedure $\delta$ is  stable with respect to $\Theta$. Assume also  that    $n/p$, $\|\beta\|^2$, $\frac{{{{\left\| \theta_f  \right\|}^2}}}{{\var[f(X)]}}$
 and
$\frac{{E\left( {{{\left\| b \right\|}^2}} \right)}}{{n\var[f(X)]}}$ are bounded, and
$\frac{{\left\| {\bf{A}} \right\|_F^2}}{{{n^2}}} \to 0$, $\frac{{\left\| {\bf{B}} \right\|_F^2}}{{{n^2\var[f(X)]}}}\to~0$, 
 and $\frac{{\left\| {\bf{C}} \right\|_F^2}}{{\{n\var[f(X)]\}^2}} \to 0.$
 In addition, assume that the first four moments of $ T_h$, $T_{\hat h}$, $T_f$, and $T_{\hat f}$ are bounded.
Then, 
\begin{equation}
    \sqrt n \left[    T_{\hat h}  -  T_h   \right]\overset{p}{\rightarrow}~0.
\end{equation}
\end{proposition}
Our proof of  Proposition \ref{prop:singel_asymptotic} shows a slightly stronger result: the proposed estimator $T_{\hat h}$  is also asymptotically equivalent to $T_f$, the oracle-estimator that originally knows the set of indices $\Theta$.

We now discuss the assumptions of Proposition \ref{prop:singel_asymptotic}.
In Remark~\ref{remark: forb_A}, several sufficient conditions implying that $\frac{{\left\| {\bf{A}} \right\|_F^2}}{{{n^2}}} \to 0$  were presented . Similarly, in Remark \ref{remark:E_norm2_b} we show that if the  covariates $X_{ij}$, for $j=1,...,p$, $i=1,...,n,$ and the response $Y$  are bounded, then so is
 $\frac{{E\left( {{{\left\| b \right\|}^2}} \right)}}{{n\var[f(X)]}}$. It is also shown, that if in addition $|\Theta|$ is bounded and $\var[f(X)]$ is bounded away from zero, then 
$\frac{{\left\| {\bf{B}} \right\|_F^2}}{{{n^2\var[f(X)]}}} \to 0$ and  $\frac{{\left\| {\bf{C}} \right\|_F^2}}{{{\{n\var[f(X)]\}^2}}} \to 0.$
   Proposition 4 in \citet{livne2021improved}  shows that   
$\frac{{\left\| {\bf{B}} \right\|_F^2}}{{{n^2\var[f(X)]}}} \to 0$ and  $\frac{{\left\| {\bf{C}}  \right\|_F^2}}{{{\{n\var[f(X)]\}^2}}} \to 0$
hold also when $|\Theta|$ is unbounded, but with additional conditions on linearity and independence of the covariates. It is also shown there that under those assumptions, $\frac{{{{\left\| \theta_f  \right\|}^2}}}{{\var[f(X)]}}$ is bounded.
In simulations, which are not presented here, we observed that these conditions also hold for various non-linear models.

Notice that the zero-estimator $Z_g$, which is not based on any covariate-selection procedure, is just a special case of $Z_h$ when $\delta$ selects all the $p$ covariates, i.e., ${\bf{S}}_\delta = \{1,\dots p \}$.   
Hence, if the conditions of Proposition~\ref{prop:singel_asymptotic} hold, then  $ \sqrt n \left[    T_{\hat g}  -  T_g   \right]\overset{p}{\rightarrow}~0,$
where $T_{\hat g}$ and $T_g$ are given in
\eqref{eq:T_g_hat} and \eqref{eq:T_g}, respectively.

\section{Simulations Results} \label{sim_res}

 In this section, we   illustrate the performance of the proposed estimators using simulations. Specifically,   we compare   the naive estimator $\hat\tau^2$  and  the improved estimators $T_{\hat g}$ and $T_{\hat h}$ which are defined  in \eqref{estimates}, \eqref{eq:T_g_hat}, and   Algorithm \ref{alg:1}, respectively. The code for reproducing the results of this section and the next section
(\ref{gener_es}) 
  is available at \url{https://t.ly/dwJg}


For demonstration purposes, we consider a  setting in which  $K$ entries of the vector $\beta$ are
relatively large (in absolute value), and all other entries are   small.
The proportion of the signal in those $K$ entries is defined  as the \emph{sparsity level}  of the vector $\beta$. Next,   we study different sparsity levels by defining the following non-linear model,
  \begin{equation}\label{non_linear_model}
{Y_i} = {\gamma _L}\sum\limits_{j \in \Theta} {\left[ {{X_{ij}} + \sin \left( {{X_{ij}}} \right)} \right]}  + {\gamma _S}\sum\limits_{j \notin \Theta } {\left[ {{X_{ij}} + \sin \left( {{X_{ij}}} \right)} \right]}  + {\xi_i}, \quad\quad  i=1,\dots,n,   
  \end{equation}
where
\({\gamma _L} \equiv {\left\{ {\frac{{\eta {\tau ^2}}}{{k{{\left( {1 + E\left[ {X\sin \left( X \right)} \right]} \right)}^2}}}} \right\}^{1/2}}\) 
, 
\({\gamma _S} \equiv {\left\{ {\frac{{{\tau ^2}\left( {1 - \eta } \right)}}{{\left( {p - k} \right){{\left( {1 + E\left[ {X\sin \left( X \right)} \right]} \right)}^2}}}} \right\}^{1/2}}\), and $\Theta$ is the set of the largest~$K$ entries of the vector $\beta.$ 
The model has two parameters,   $\tau^2$ and $\eta$, that vary across the different simulation scenarios.
The covariates were generated  from the centered exponential distribution, i.e.,   $X_{ij}\overset{iid}{\sim} \text{Exp}(1)-1,$  $i=1,\dots,n$,  $j=1,\dots, p$. The noise level $\xi_i$ was generated from the standard normal distribution. 
 One can verify that under the above model   $\beta _j^2 = \frac{{\eta {\tau ^2}}}{K}$ for $j \in \Theta$, and that
$\beta _j^2 = \frac{{{\tau ^2}\left( {1 - \eta } \right)}}{{\left( {p - K} \right)}}\,$ for \(j \notin \Theta.\)
   Define \(\tau _\Theta^2 \equiv \sum\limits_{j \in \Theta}^{} {\beta _j^2} \). From the above definitions, it follows that~\(\eta=~\tau _\Theta^2/{\tau ^2}\).
    The parameter $\eta$ is the proportion of  signal that is captured by the set $\Theta,$ which is the sparsity level as defined above. The case of  full sparsity, where the entire   signal level $\tau^2$ comes only from the  set $\Theta$, corresponds to $\eta=1$, and  is not assumed here.

 We fix $n=p=300$ and $K=6$.
For each combination of the parameters  $\tau^2\in \{1,2\}$ and \(\eta \in\{ {0.1, 0.3, 0.5, 0.7, 0.9} \},\) we generated  100  independent datasets from model \eqref{non_linear_model} and estimated $\tau^2$  using the different estimators.  
The covariate-selection procedure $\delta$ that was used in the estimator $T_{\hat h}$ is defined in Remark \ref{selection_algorithm} in the Appendix.
  
Figure \ref{figure1} plots the RMSE of each estimator as a function of the sparsity  level $\eta$ and the signal level~$\tau^2.$ It is demonstrated that
 the estimators $T_{\hat g}$ and $T_{\hat h}$ improve (i.e., lower or equal RMSE) the naive estimator 
  in all settings.
 The improved estimators are complementary to each other, i.e.,
for small values of $\eta$ the estimator $T_{\hat g}$  performs better than   $T_{\hat h},$ and the opposite occurs for large values of $\eta.$ 
This is expected since 
 when the sparsity level $\eta$ is small, 
 the improvement of $T_{\hat h}$ is smaller as it ignores much of the signal that lies outside of the set $\Theta.$ 
On the other hand, when a
large portion of the signal $\tau^2$ is captured by only the few  covariates in $\Theta$, it is sufficient to make use of only these covariates in the  zero-estimator term, and the improvement of $T_{\hat h}$ is greater.

Table \ref{table1} shows the RMSE, bias, standard error, and the relative improvement, for the different estimators. 
  It can be observed that the degree of improvements depends on  the sparsity level of the data $\eta$.
 For example, when $\tau^2=1$ and sparsity level is low ($\eta=0.1$),  the estimator  $T_{\hat g}$ improves the naive estimator by   $11\%$, while the estimator $T_{\hat h}$  presents a  similar performance to the naive estimator.
 On the other hand, when the sparsity level is high ($\eta=0.9$), the estimator $T_{\hat h}$ improves the naive by $11\%$, while $T_{\hat g}$ presents a  similar performance to the naive estimator, as expected.
 Notice that  when $\tau^2=2$ these improvements are even more substantial.

\begin{figure}[H]
  \centering
 \includegraphics[width=0.8\textwidth]{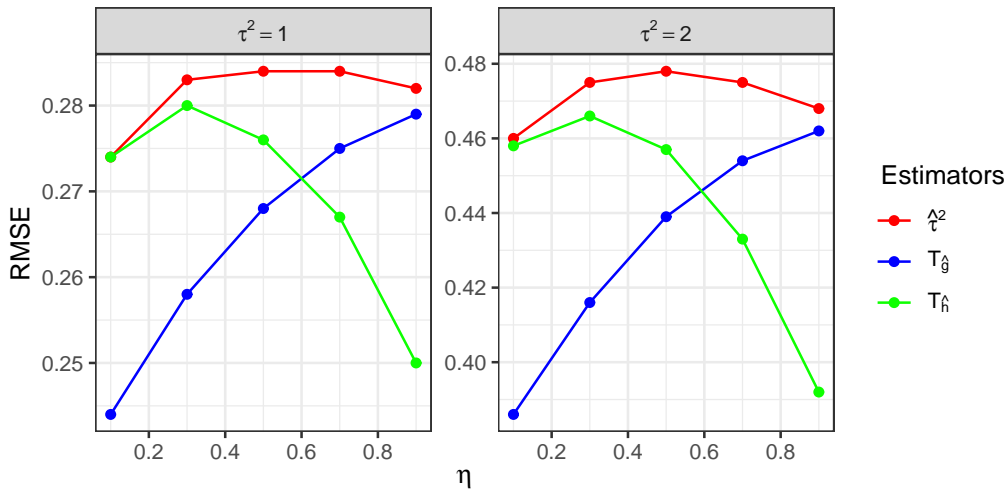}
 \captionsetup{font=footnotesize}
\caption{
Root mean square error (RMSE) for the proposed estimators.   The x-axis stands for the sparsity level~$\eta$. 
}
\label{figure1}
\end{figure}

\newpage

\begin{table}[H] 
\small
\captionsetup{font=scriptsize}
\caption{Summary statistics for the proposed estimators;  $n=p=300.$ Bias, standard error (SE),  root mean square error (RMSE) and percentage change from the naive estimator (in terms of RMSE) are shown. The table results were computed over $100$ simulated  datasets for each setting. 
An estimate for the standard deviation of RMSE ($\hat\sigma_{RMSE}$) was calculated using the delta method.}
 \label{table1} \par
\resizebox{\linewidth}{!}{
 \renewcommand{\arraystretch}{0.8}   
\resizebox{1.5cm}{!}{
\begin{tabular}{|cccccccc|} \hline 

$\eta$ & $\tau^2$  &   Estimator &  Bias & SE & RMSE & \% Change & $\hat\sigma_{RMSE}$   \\ \hline
10\% & 1   & $\hat\tau^2$ & -0.03 & 0.274 & 0.274 & 0.00 & 0.018 \\ 
10\% & 1   & $T_{\hat g}$ & 0.04 & 0.242 & 0.244 & -10.95 & 0.016 \\ 
10\% & 1   & $T_{\hat h}$ & -0.03 & 0.274 & 0.274 & 0.00 & 0.017 \\
\hline
30\% & 1   & $\hat\tau^2$ & -0.02 & 0.284 & 0.283 & 0.00 & 0.019 \\ 
30\% & 1   & $T_{\hat g}$ & 0.05 & 0.255 & 0.258 & -8.83 & 0.018 \\ 
30\% & 1   & $T_{\hat h}$ & 0.00 & 0.282 & 0.28 & -1.06 & 0.017 \\ 
\hline
50\% & 1   & $\hat\tau^2$ & 0.00 & 0.286 & 0.284 & 0.00 & 0.021 \\ 
50\% & 1   & $T_{\hat g}$ & 0.05 & 0.264 & 0.268 & -5.63 & 0.020 \\ 
50\% & 1   & $T_{\hat h}$ & 0.02 & 0.277 & 0.276 & -2.82 & 0.017 \\ 
\hline
70\% & 1   & $\hat\tau^2$ & 0.01 & 0.285 & 0.284 & 0.00 & 0.022 \\ 
70\% & 1   & $T_{\hat g}$ & 0.05 & 0.272 & 0.275 & -3.17 & 0.021 \\ 
70\% & 1   & $T_{\hat h}$ & 0.04 & 0.265 & 0.267 & -5.99 & 0.017 \\ 
\hline
90\% & 1   & $\hat\tau^2$ & 0.03 & 0.281 & 0.282 & 0.00 & 0.021 \\ 
90\% & 1   & $T_{\hat g}$ & 0.05 & 0.276 & 0.279 & -1.06 & 0.020 \\ 
90\% & 1   & $T_{\hat h}$ & 0.07 & 0.242 & 0.25 & -11.35 & 0.015 \\ 
\hline\hline
10\% & 2   & $\hat\tau^2$ & -0.06 & 0.458 & 0.46 & 0.00 & 0.030 \\ 
10\% & 2   & $T_{\hat g}$ & 0.08 & 0.379 & 0.386 & -16.09 & 0.024 \\ 
10\% & 2   & $T_{\hat h}$ & -0.05 & 0.457 & 0.458 & -0.43 & 0.029 \\
\hline
30\% & 2   & $\hat\tau^2$ & -0.03 & 0.476 & 0.475 & 0.00 & 0.033 \\ 
30\% & 2   & $T_{\hat g}$ & 0.09 & 0.408 & 0.416 & -12.42 & 0.029 \\ 
30\% & 2   & $T_{\hat h}$ & -0.01 & 0.469 & 0.466 & -1.89 & 0.029 \\ 
\hline
50\% & 2   & $\hat\tau^2$ & -0.01 & 0.481 & 0.478 & 0.00 & 0.036 \\ 
50\% & 2   & $T_{\hat g}$ & 0.09 & 0.431 & 0.439 & -8.16 & 0.033 \\ 
50\% & 2   & $T_{\hat h}$ & 0.03 & 0.458 & 0.457 & -4.39 & 0.028 \\ 
\hline
70\% & 2   & $\hat\tau^2$ & 0.02 & 0.477 & 0.475 & 0.00 & 0.038 \\ 
70\% & 2   & $T_{\hat g}$ & 0.09 & 0.448 & 0.454 & -4.42 & 0.035 \\ 
70\% & 2   & $T_{\hat h}$ & 0.08 & 0.429 & 0.433 & -8.84 & 0.026 \\ 
\hline
90\% & 2   & $\hat\tau^2$ & 0.05 & 0.468 & 0.468 & 0.00 & 0.035 \\ 
90\% & 2   & $T_{\hat g}$ & 0.08 & 0.456 & 0.462 & -1.28 & 0.034 \\ 
90\% & 2   & $T_{\hat h}$ & 0.12 & 0.376 & 0.392 & -16.24 & 0.023 \\ 
\hline

\end{tabular}}
}
\end{table}

 \section{ Generalization to Other Estimators
 }\label{gener_es}
The suggested methodology in this paper is not limited  to improving only the naive estimator, but can also be generalized to other  estimators. As before, the key idea is to use a zero-estimator that is correlated with an initial estimator of $\tau^2$ in order to reduce its variance.   Unlike the naive estimator~$\hat\tau^2$, which has by a closed-form expression, other  estimators, such as the EigenPrism estimator \citep{janson2017eigenprism}, are computed numerically by solving a convex optimization problem. For a given zero-estimator, this makes the task of estimating the  optimal-coefficient~$c^*$  more challenging than before.
To overcome this challenge, we  approximate the optimal coefficient $c^*$ using bootstrap samples.
This is described in the following algorithm.

\vspace{10mm}

\begin{algorithm}[H]\label{alg_emp}
\SetAlgoLined

\vspace{0.4 cm}
\textbf{Input:} 
 A dataset \(\left( {{{\bf{X}}_{n \times p}},{Y_{n \times 1}}} \right)\), an estimation procedure  $\tilde{\tau}^2,$   and a covariate-selection procedure~$\delta$.
\begin{enumerate}
 \item Apply the procedure $\delta$ to the dataset \(\left( {{{\bf{X}}_{n \times p}},{Y_{n \times 1}}} \right)\) to obtain  \({{\bf{S}}_{\delta}}.\) 
\item Apply the procedure $\tilde{\tau}^2$  to the dataset. 

  \item   Calculate the zero-estimator
 \({Z_h} = \frac{1}{n}\sum\limits_{i = 1}^n {h\left( {{X_i}} \right)} \),
     where   \(h\left( {{X_i}} \right) = \sum\limits_{j < j' \in {{\bf{S}}_\delta }}^{} {{X_{ij}}{X_{ij'}}}.\)
    
     \item \textbf{Bootstrap step:}  
      \begin{itemize}
        \item  Sample    $n$ observations at random  from \(\left( {{{\bf{X}}_{n \times p}},{Y_{n \times 1}}} \right)\), with replacement, to obtain a bootstrap dataset.
        \item Repeat steps 2 and 3   based on the bootstrap dataset. 
    \end{itemize}
    The bootstrap step is repeated $M$ times  in order to produce
    $(\tilde{\tau}^2)^{*1},...,(\tilde{\tau}^2)^{*M}$ and 
    \(Z_h^{*1},...,Z_h^{*M}.\)
    \item   Approximate the  coefficient  	 
    $\tilde{c}_h^* =  \frac{{\widehat {\cov\left( {\tilde{\tau}^2,{Z_h}} \right)}}}{{\var\left( {{Z_h}} \right)}}$
    where \(\widehat {{\cov} \left(  \cdot  \right)}\) denotes the empirical covariance from the bootstrap samples. 
  \end{enumerate}
\KwResult{Return the empirical estimator
$T_{\tilde h}\equiv\tilde \tau^2 - \tilde c_h^* Z_h.$
\vspace{0.4 cm}
}
  \caption{Empirical Estimators}
\end{algorithm}
\vspace{10mm}

In the special case when $\delta$ selects all the $p$ covariates, i.e., ${\bf{S}}_\delta = \{1,\dots p \}$,  we use the notations $Z_g$ and $\tilde c_g^*$ rather than $Z_h$ and $\tilde c_h^*$, respectively, i.e., $T_{\tilde g}\equiv\tilde \tau^2 - \tilde c_g^* Z_g.$

We illustrate the improvement obtained by Algorithm~\ref{alg_emp} by choosing $\tilde\tau^2$ to be  the EigenPrism   procedure \citep{janson2017eigenprism}, but other estimators can be used as well.
We consider the same  setting as in Section \ref{sim_res}.  The number of bootstrap samples is $M= 100.$ 

\begin{figure}[H]
  \centering
   \includegraphics[width=0.8\textwidth]{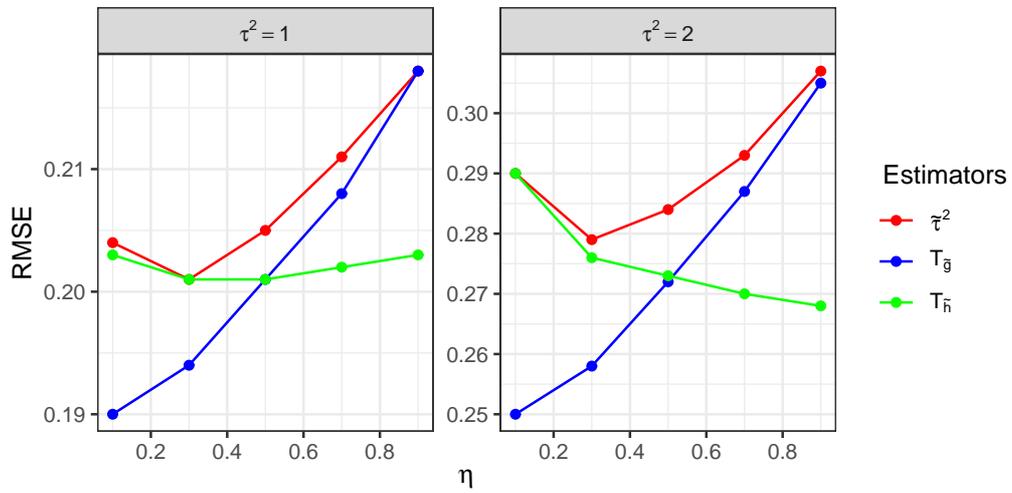}
   \captionsetup{font=footnotesize}
\caption{
Root mean square error (RMSE) for the proposed estimators.   The x-axis stands for the sparsity level~$\eta$.
}
\label{figure_eigen}
\end{figure}

The  simulation results appear in Table \ref{table_emp} and Figure \ref{figure_eigen}.
 Both estimators $T_{\tilde h}$ and $T_{\tilde g}$   show an improvement over the EigenPrism estimator $\tilde\tau^2.$ 
 The results here are fairly similar to the results  shown for the naive estimator in Section \ref{sim_res}, with just a smaller degree of improvement. As before, 
 the improved estimators $T_{\tilde h}$ and $T_{\tilde g}$ are complementary to each other, i.e.,
for small values of $\eta$ the estimator~$T_{\tilde g}$  performs better than   $T_{\tilde h},$ and the opposite occurs for large values of $\eta.$

\vspace{10mm}

\begin{table}[H] 
\small
\captionsetup{font=scriptsize}
\caption{Summary statistics for the EigenPrism-related estimators;  $n=p=300.$ Biases, standard errors (SE) and the root mean square errors (RMSE)  of the different estimators, computed over $100$ independent datasets for each setting. The relative improvement over the EigenPrism estimator $\tilde\tau^2$ (in terms of RMSE) is also shown.      
An estimate for the standard deviation of RMSE ($\hat\sigma_{RMSE}$) was calculated using the delta method.}
  \label{table_emp}  \par
  
  \resizebox{\linewidth}{!}{
 \renewcommand{\arraystretch}{0.8}   
\resizebox{1.5cm}{!}{
\begin{tabular}{|cccccccc|} \hline 
$\eta$ & $\tau^2$  &   Estimator &  Bias & SE & RMSE & \% Change & $\hat\sigma_{RMSE}$   \\ \hline
10\% & 1  & $\tilde\tau^2$ & 0.01 & 0.204 & 0.204 & 0.00 & 0.014 \\ 
10\% & 1  & $T_{\tilde g}$ & 0.01 & 0.19 & 0.19 & -6.86 & 0.012 \\ 
10\% & 1  & $T_{\tilde h}$ & 0.01 & 0.204 & 0.203 & -0.49 & 0.014 \\ 
\hline 
30\% & 1  & $\tilde\tau^2$ & 0.01 & 0.202 & 0.201 & 0.00 & 0.014 \\ 
30\% & 1  & $T_{\tilde g}$ & 0.01 & 0.195 & 0.194 & -3.48 & 0.014 \\ 
30\% & 1  &$T_{\tilde h}$ & 0.01 & 0.201 & 0.201 & 0.00 & 0.014 \\ 
\hline 
50\% & 1  & $\tilde\tau^2$ & 0.00 & 0.206 & 0.205 & 0.00 & 0.014 \\ 
50\% & 1  & $T_{\tilde g}$ & 0.00 & 0.202 & 0.201 & -1.95 & 0.014 \\ 
50\% & 1  & $T_{\tilde h}$ & 0.01 & 0.202 & 0.201 & -1.95 & 0.015 \\ 
\hline 
70\% & 1  & $\tilde\tau^2$ & 0.00 & 0.212 & 0.211 & 0.00 & 0.015 \\ 
70\% & 1  & $T_{\tilde g}$ & 0.00 & 0.209 & 0.208 & -1.42 & 0.015 \\ 
70\% & 1  & $T_{\tilde h}$ & 0.00 & 0.203 & 0.202 & -4.27 & 0.015 \\ 
\hline 
90\% & 1  & $\tilde\tau^2$ & -0.01 & 0.219 & 0.218 & 0.00 & 0.016 \\ 
90\% & 1  & $T_{\tilde g}$ & -0.01 & 0.218 & 0.218 & 0.00 & 0.016 \\ 
90\% & 1  & $T_{\tilde h}$ & 0.00 & 0.204 & 0.203 & -6.88 & 0.015 \\ 

\hline\hline 
10\% & 2  & $\tilde\tau^2$ & 0.00 & 0.291 & 0.29 & 0.00 & 0.020 \\ 
10\% & 2  & $T_{\tilde g}$ & 0.02 & 0.251 & 0.25 & -13.79 & 0.016 \\ 
10\% & 2  & $T_{\tilde h}$ & 0.01 & 0.291 & 0.29 & 0.00 & 0.020 \\ 
\hline 
30\% & 2  & $\tilde\tau^2$ & 0.03 & 0.279 & 0.279 & 0.00 & 0.019 \\ 
30\% & 2  & $T_{\tilde g}$ & 0.03 & 0.258 & 0.258 & -7.53 & 0.018 \\ 
30\% & 2  & $T_{\tilde h}$ & 0.03 & 0.276 & 0.276 & -1.08 & 0.019 \\ 
\hline 
50\% & 2  & $\tilde\tau^2$ & 0.02 & 0.285 & 0.284 & 0.00 & 0.019 \\ 
50\% & 2  & $T_{\tilde g}$ & 0.01 & 0.273 & 0.272 & -4.23 & 0.019 \\ 
50\% & 2  & $T_{\tilde h}$ & 0.02 & 0.274 & 0.273 & -3.87 & 0.020 \\ 

\hline 
70\% & 2  & $\tilde\tau^2$ & 0.01 & 0.294 & 0.293 & 0.00 & 0.020 \\ 
70\% & 2  & $T_{\tilde g}$ & 0.00 & 0.289 & 0.287 & -2.05 & 0.019 \\ 
70\% & 2  & $T_{\tilde h}$ & 0.01 & 0.271 & 0.27 & -7.85 & 0.019 \\ 
\hline 
90\% & 2  & $\tilde\tau^2$ & -0.01 & 0.308 & 0.307 & 0.00 & 0.022 \\ 
90\% & 2  & $T_{\tilde g}$ & -0.01 & 0.307 & 0.305 & -0.65 & 0.022 \\ 
90\% & 2  & $T_{\tilde h}$ & 0.00 & 0.269 & 0.268 & -12.7 & 0.019 \\ 
\hline  
\end{tabular}}
} 
\end{table}

   \newpage  
\section{Discussion and Future Work }\label{discuss}
In this work, we proposed  a zero-estimator approach for improving  estimation of the signal and noise levels explained by a set of covariates in a   high-dimensional regression setting when the covariate distribution is known.  
We  presented  theoretical properties of the naive estimator $\hat\tau^2,$  and the proposed improved estimators  $T_{\hat h}$ and $T_{\hat g}$.
In a simulation study, we  demonstrated that the zero-estimator approach  leads to a significant reduction in the RMSE. 
Our method does not  rely on sparsity assumptions of the regression coefficient vector, normality of the covariates,  or linearity of $E(Y|X)$. The goal in this work is to estimate the signal coming from the best linear function of the covariates, which is a model-free quantity.
Our simulations demonstrate that our approach can be generalized to improve other estimators as well.

 We suggest the following directions for future work. 
 One natural extension is to relax the  assumption of  known covariate distribution to allow for a more general setting.
 This may be studied under the semi-supervised setting where one has access to a large amount of unlabeled data  ($N \gg n$) in order to obtain theoretical results as a function of $N.$  Another possible future research
 might be  to extend the proposed approach to  generalized linear models (GLM)  such as logistic and Poisson regression, or survival models.

\bibhang=1.7pc
\bibsep=2pt
\fontsize{9}{14pt plus.8pt minus .6pt}\selectfont
\renewcommand\bibname{\large \bf References}
\expandafter\ifx\csname
natexlab\endcsname\relax\def\natexlab#1{#1}\fi
\expandafter\ifx\csname url\endcsname\relax
  \def\url#1{\texttt{#1}}\fi
\expandafter\ifx\csname urlprefix\endcsname\relax\def\urlprefix{URL}\fi

  \bibliographystyle{chicago}      
  \bibliography{bib_ref}   








\newpage

\section{Appendix}\label{appnd}

 \begin{proof}[Proof of Proposition~\ref{prop:var_naive_tau2}]$ $\newline
Let ${W_i} = {\left( {{W_{i1}},...,{W_{ip}}} \right)^T}$ and notice that  $\hat\tau^2= \frac{1}{{n\left( {n - 1} \right)}} \sum\limits_{{i_1} \ne {i_2}}^n\sum\limits_{j = 1}^p {{W_{{i_1}j}}{W_{{i_2}j}}} $ is a U-statistic of order~2 with the kernel \(h_{\tau}\left( {{W_1},{W_2}} \right) = W_1^T{W_2} = \sum\limits_{j = 1}^p {{W_{1j}}{W_{2j}}}.\)
 
 By Theorem 12.3 in \citet{van2000asymptotic},
 \begin{equation}\label{naive_var_form}
  {\var} \left( {{{\hat \tau }^2}} \right) = \frac{{4\left( {n - 2} \right)}}{{n\left( {n - 1} \right)}}{\zeta _1} + \frac{2}{{n\left( {n - 1} \right)}}{\zeta _2},   
  \end{equation}  where \({\zeta _1} \equiv {\cov} \left[ {h\left( {{W_1},{W_2}} \right),h\left( {{W_1},{{\widetilde{W}_2}}} \right)} \right]\)
  and 
  ${\zeta _2} \equiv {\var} \left[ {h\left( {{W_1},{W_2}} \right)} \right],$
   and where $\widetilde{W}_2$ is an independent copy of $W_2$.
 Define the $p \times p$ matrix  
   \({\bf{A}} = E\left( {{{W}}{W}^T} \right)\)   and   notice that
  \begin{align*}
  {\zeta _1} &\equiv \cov\left[ {h\left( {{{W}_1},{{W}_2}} \right),h\left( {{{W}_1},{{\widetilde {W}}_2}} \right)} \right]
  \\ &= \sum\limits_{j,j'}^p {\cov\left( {{W_{1j}}{W_{2j}},{W_{1j'}}{{\widetilde W}_{2j'}}} \right)}  = \sum\limits_{j,j'}^p {\left( {{\beta _j}{\beta _{j'}}E\left[ {{W_{1j}}{W_{1j'}}} \right] - \beta _j^2\beta _{j'}^2} \right)}\\
  &= {\beta ^T}{\bf{A}}\beta  - {\left\| \beta  \right\|^4},    
  \end{align*}
and 
\begin{align*}
{\zeta _2} &\equiv {\cov} \left[ {h\left( {{{W}_1},{{W}_2}} \right),h\left( {{{W}_1},{{W}_2}} \right)} \right]\\
&= \sum\limits_{j,j'}^{} {{\cov} \left( {{W_{1j}}{W_{2j}},{W_{1j'}}{W_{2j'}}} \right)}
=\sum\limits_{j,j'}^{}\left( {{{\left( {E\left[ {{W_{1j}}{W_{1j'}}} \right]} \right)}^2} - \beta _j^2\beta _{j'}^2}\right)\\  
&= \left\| {\bf{A}} \right\|_F^2 - {\left\| \beta  \right\|^4},    
\end{align*}
where  $\|\textbf{A}\|_F^{2}$ is the Frobenius norm of $\textbf{A}.$ Thus, \eqref{eq:var_naive} follows from  (\ref{naive_var_form}). 
\end{proof}

\begin{proof}[Proof of Proposition \ref{consistency_naive}]$ $\newline
Since $\hat\tau^2$ is an unbiased estimator of $\tau^2,$ then it is enough to prove that  $\var[\hat\tau^2]\xrightarrow{n\rightarrow\infty}0$.
By \eqref{eq:var_naive} it is enough to require that $\frac{\beta^T\textbf{A}\beta}{n}\xrightarrow{n\rightarrow\infty}0$ and $\frac{\|\textbf{A}\|^{2}_F}{n^2}\xrightarrow{n\rightarrow\infty}0$. The latter is assumed  and we now show that the former also holds true.\\
Let $\lambda_{1{\bf A}}\geq...\geq\lambda_{p{\bf A}}$ be the eigenvalues of $\textbf{A}$  
and notice that $\textbf{A}$ is symmetric.
Thus,
$$\frac{1}{n^2}(\lambda_{1{\bf A}})^2 \leq \frac{1}{n^2}\sum_{j=1}^{p}\lambda_{j{\bf A}}^2=\frac{1}{n^2} tr(\textbf{A}^2)=\frac{1}{n^2}\|\textbf{A}\|_F^2.$$
Since we assume that $\frac{\|\textbf{A}\|^{2}_F}{n^2}\xrightarrow{n\rightarrow\infty}0$, we can conclude that $\frac{\lambda_{1{\bf A}}}{n}\xrightarrow{n\rightarrow\infty}0$. 
Now, recall that the maximum of the quadratic form ${a^T}{\bf{A}}a$  satisfies
${\lambda_{1{\bf A}}} = \mathop {\max }\limits_a \,\frac{{{a^T}{\bf{A}}a}}{{{{\left\| a \right\|}^2}}}.$ Hence, 
$$\frac{1}{n}\beta^T\textbf{A}\beta\equiv\frac{1}{n}\|\beta\|^2[(\frac{\beta}{\|\beta\|})^T\textbf{A}\frac{\beta}{\|\beta\|}]\leq\frac{1}{n}\|\beta\|^2\lambda_1 = \frac{\lambda_1}{n}\tau^2\xrightarrow{n\rightarrow\infty}0,
$$
where the last limit follows from the assumption that $\tau^2=O(1),$ and from the fact that $\frac{\lambda_1}{n}\xrightarrow{n\rightarrow\infty}0.$\\
We conclude that $\var[\hat\tau^2] \xrightarrow{n\rightarrow\infty} 0$. 
\end{proof}

\begin{proof}[Proof of Proposition \ref{proposition: consistent_sigma2}]$ $\newline
Recall that  \eqref{eq: sigma2_estimator} states 
 that ${\hat \sigma ^2} \equiv \hat \sigma _Y^2 - {\hat \tau ^2}$.
Thus,
\begin{align}\label{eq:var sigma^2}
\var\left( {{{\hat \sigma }^2}} \right) = \var\left( {{{\hat \tau }^2}} \right) + \var\left( {\hat \sigma _Y^2} \right)   - 2\cov\left(\hat \sigma _Y^2 ,{\hat \tau }^2\right).    \end{align}
The variance of $\hat\tau^2$  is given in \eqref{eq:var_naive}. By standard U-statistic calculations (see e.g., Example 1.8 in \citet{bose2018u}), the variance of $\hat\sigma_Y^2$ is 
\begin{align*}
\text{Var}(\hat\sigma_{Y}^2)
=\frac{4(n-2)}{n(n-1)}\psi_1+\frac{2}{n(n-1)}\psi_2,
\end{align*}
where
$\psi_1 \equiv\frac{\mu_4-\sigma_Y^4}{4}$ and $\psi_2 \equiv\frac{\mu_4+\sigma_Y^4}{2}.$ This can be further simplified to obtain
\begin{align}\label{eq:var_var_y}
\var\left( {\hat \sigma _Y^2} \right) =  \frac{1}{n}{\mu _4} - \frac{{\left( {n - 3} \right)}}{{n\left( {n - 1} \right)}}\sigma _Y^4.
\end{align}

We now calculate the covariance between  and $\hat\sigma_Y^2$ and $\hat\tau^2.$
 Write,
 \begin{align*}
 \cov\left( {\hat \sigma _Y^2,{{\hat \tau }^2}} \right) = \cov\left[ {\frac{1}{{n\left( {n - 1} \right)}}\sum\limits_{{i_1} \ne {i_2}}^{} {{{\left( {{Y_{{i_1}}} - {Y_{{i_2}}}} \right)}^2}/2,\frac{1}{{n\left( {n - 1} \right)}}\sum\limits_{j = 1}^p {\sum\limits_{{i_1} \ne {i_2}}^{} {{W_{{i_1}j}}{W_{{i_2}j}}} } } } \right]\\ = \frac{1}{{2{n^2}{{\left( {n - 1} \right)}^2}}}\sum\limits_{j = 1}^p {\sum\limits_{{i_1} \ne {i_2}}^{} {\sum\limits_{{i_3} \ne {i_4}}^{} {\cov\left[ {{W_{{i_3}j}}{W_{{i_4}j}}},{{\left( {{Y_{{i_1}}} - {Y_{{i_2}}}} \right)}^2} \right]} }}. 
 \end{align*}
 The covariance above is different from zero either when one of $\{i_1, i_2 \}$ is equal to one of $\{i_3,i_4 \}$, or when  $\{i_1, i_2 \}$ is equal to $\{i_3,i_4 \}$. There are $4n(n-1)(n-2)$ quadruples $(i_1,i_2,i_3,i_4)$ for the first case, and $2n(n-1)$ for the second case. Therefore,
\begin{align}\label{eq:cov_var_y_naive}
\cov[\hat\sigma_{Y}^2,\hat{\tau}^2] = \frac{2(n-2)}{n(n-1)}\sum_{j=1}^{p}\cov[W_{1j}{W_{2j}},(Y-Y_3)^2]+\frac{1}{n(n-1)}\sum_{j=1}^{p}\cov[W_{1j}W_{2j},(Y_1-{Y_2})^2].
\end{align}
The first covariance term in (\ref{eq:cov_var_y_naive}) is
\begin{align*}
&\cov[W_{1j}W_{2j},(Y_1-Y_3)^2]  = E\Big( [W_{1j}W_{2j}-E(W_{1j})E(W_{2j})][(Y_1-Y_3)^2-2\sigma_Y^2]  \Big)\\
&=E\Big( [W_{1j}W_{2j}-\beta_j^2][(Y_1-\alpha-(Y_3-\alpha))^2-2\sigma_Y^2]  \Big)\\
&=E\Big( [W_{1j}W_{2j}-\beta_j^2][\{(Y_1-\alpha)^2-\sigma_Y^2\}-2(Y_1-\alpha)(Y_3-\alpha)+(Y_3-\alpha)^2-\sigma_Y^2]  \Big)\\
&=E\Big( [W_{1j}W_{2j}-\beta_j^2][\{(Y_1-\alpha)^2-\sigma_Y^2\}]\Big)\\
&=E\left[ {{W_{1j}}{{\left( {Y_1 - \alpha } \right)}^2}} \right]\beta _j^2 - \beta _j^2\sigma _Y^2.
\end{align*}
The second covariance term of  (\ref{eq:cov_var_y_naive}) is
\begin{align*}
&\cov\left[ {W_{1j}W_{2j},{{\left( Y_1-Y_2 \right)}^2}} \right] = E\left\{ {\left[ {W_{1j}W_{2j} - E\left( {W_{1j}W_{2j}} \right)} \right]\left[ {{{\left( Y_1-Y_2 \right)}^2} - E\left\{ {{{\left( Y_1-Y_2 \right)}^2}} \right\}} \right]} \right\}\\ 
&= E\left\{ {\left( {W_{1j}W_{2j} - \beta _j^2} \right)\left[ {{{\left( Y_1-Y_2 \right)}^2} - 2\sigma _Y^2} \right]} \right\}\\
&= E\left\{ {\left( {W_{1j}W_{2j} - \beta _j^2} \right)\left[ {{{\left\{ {\left( {Y_1 - \alpha } \right) - \left( {Y_2 - \alpha } \right)} \right\}}^2} - 2\sigma _Y^2} \right]} \right\}\\
&= E\left\{ {\left( {W_{1j}W_{2j} - \beta _j^2} \right)\left[ {\left\{ {{{\left( {Y_1 - \alpha } \right)}^2} - \sigma _Y^2} \right\} - 2\left( {Y_1 - \alpha } \right)\left( {Y_2 - \alpha } \right) + \left\{ {{{\left( {Y_2 - \alpha } \right)}^2} - \sigma _Y^2} \right\}} \right]} \right\}\\
&= 2E\left\{ {\left[ {{{\left( {Y_1 - \alpha } \right)}^2} - \sigma _Y^2} \right]\left( {W_{1j}W_{2j} - \beta _j^2} \right)} \right\} - 2E\left[ {\left( {W_{1j}W_{2j} - \beta _j^2} \right)\left( {Y_1 - \alpha } \right)\left( {Y_2 - \alpha } \right)} \right]\\ 
&= 2\left\{ {E\left[ {{{\left( {Y_1 - \alpha } \right)}^2}{W_{1j}}} \right]{\beta _j} - \sigma _Y^2\beta _j^2} \right\} - 2E\left[ {{W_{1j}}\left( {Y_1 - \alpha } \right) W_{2j}\left( {Y_2 - \alpha } \right)} \right]\\ 
&= 2E\left[ {{{\left( {Y_1 - \alpha } \right)}^2}{W_{1j}}} \right]{\beta _j} - 2\sigma _Y^2\beta _j^2 - 2{\left\{ {E\left[ {{W_{1j}}\left( {Y_1 - \alpha } \right)} \right]} \right\}^2},
 \end{align*}
Hence, by \eqref{eq:cov_var_y_naive} we have
\begin{align*}
\cov\left( {\hat \sigma _Y^2,{{\hat \tau }^2}} \right) &= \frac{{2\left( {n - 2} \right)}}{{n\left( {n - 1} \right)}}\sum\limits_{j = 1}^p {\left\{ {E\left[ {W_{1j}{{\left( {Y_1 - \alpha } \right)}^2}} \right]{\beta _j} - \beta _j^2\sigma _Y^2} \right\}}\\
&+ \frac{1}{{n\left( {n - 1} \right)}}\sum\limits_{j = 1}^p {\left\{ {2E\left[ {{{\left( {Y_1 - \alpha } \right)}^2}W_{1j}} \right]{\beta _j} - 2\sigma _Y^2\beta _j^2 - 2{{\left\{ {E\left[ {W_{1j}\left( {Y_1 - \alpha } \right)} \right]} \right\}}^2}} \right\}},    
\end{align*}
which can be further simplified as
\begin{align}\label{eq:cov_tau_sig_Y}
\cov\left( {\hat \sigma _Y^2,{{\hat \tau }^2}} \right) = \frac{2}{n}\left( {{\pi ^T}\beta  - {\tau ^2}\sigma _Y^2} \right) - \frac{2}{{n\left( {n - 1} \right)}}\sum\limits_{j = 1}^p {{{\left\{ {E\left[ {W_{1j}\left( {Y_1 - \alpha } \right)} \right]} \right\}}^2}}, 
\end{align}
where  \(\pi  = {\left( {{\pi _1},...,{\pi _p}} \right)^T}\) and \({\pi _j} \equiv E\left[ {{{\left( {Y_1 - \alpha } \right)}^2}W_{1j}} \right]\).\\
Plugging   \eqref{eq:var_var_y} and \eqref{eq:cov_tau_sig_Y} into  \eqref{eq:var sigma^2} leads to \eqref{eq:var_sigma_2_hat}.
\end{proof}

\begin{proof}[Proof of Corollary \ref{cor:consistency_naive_sigma2}]$ $\newline
Since $\hat\sigma^2$ is an unbiased estimator of $\sigma^2,$ it is enough to prove that $\var(\hat\sigma^2)\xrightarrow{n\rightarrow\infty}~0.$ Recall that $\hat{\sigma}^2=\hat{\sigma}^2_Y-\hat{\tau}^2$. It follows that 
$\var\left( {{{\hat \sigma }^2}} \right) \le 2\var\left( {\hat \sigma _Y^2} \right) + 2\var\left( {{{\hat \tau }^2}} \right)$.
Thus, it is enough to prove that $\var(\hat\sigma_Y^2)\xrightarrow{n\rightarrow\infty}~0$ and $\var(\hat\tau^2)\xrightarrow{n\rightarrow\infty}~0.$ 
The former requires, by~\eqref{eq:var_var_y}, the assumption that $\mu_4$ is bounded and the latter holds true by  Proposition~\ref{consistency_naive}. 
\end{proof}

\begin{remark}\label{remark: forb_A}
The condition $\frac{\|{\bf A}\|_F^{2}}{n^2}\rightarrow 0$ holds in the homoskedastic linear model 
with the additional assumption that the columns of $\bf{X}$ are independent (\citet{livne2021improved}, Proposition 2). 
We now show two more examples where the condition $\frac{\|{\bf A}\|_F^{2}}{n^2}\rightarrow 0$ holds without assuming linearity.

1)   We show that if $p/n^2\rightarrow 0 $  and $Y^2 \leq C$ for some constant $C$, then $\frac{\|{\bf A}\|_F^{2}}{n^2}\rightarrow 0$.  For $a \in \mathbb{R}^p$ we have,
\begin{align}\label{eq: aAa}
{a^T}{\bf{A}}a &= {a^T}E\left( {{\bf{W}}{{\bf{W}}^T}} \right)a = E\left( {{a^T}{\bf{W}}{{\bf{W}}^T}a} \right) = E\left[ {{{\left( {{a^T}{\bf{W}}} \right)}^2}} \right]\nonumber\\ 
&= E\left[ {{{\left( Y {\sum\limits_{j = 1}^p {{a_j}{X_{ij}}} } \right)}^2}} \right] \le {C}E\left[ {{{\left( {\sum\limits_{j = 1}^p {{a_j}{X_{ij}}} } \right)}^2}} \right]\nonumber\\
&=  {C}\left\{ {\sum\limits_{j = 1}^p {a_j^2\underbrace {E\left( {X_{ij}^2} \right)}_1} } \right\} + CE\left( {\sum\limits_{j \ne j'}^{} {{a_j}{a_{j'}}} \underbrace {E\left( {{X_{ij}}{X_{j'}}} \right)}_0} \right)\nonumber\\
&= {C}{\left\| a \right\|^2},
\end{align}
where the last equality follows from ${\Sigma}={\bf I}.$
Now, let $\lambda_{1{\bf A}}\geq\lambda_2^{\bf A}\geq,...,\geq\lambda_{p{\bf A}}^{\bf A}$ be the eigenvalues of $\bf{A},$ and recall that the extrema of the quadratic form ${a^T}{\bf{A}}a$  satisfies
${\lambda_{1{\bf A}}} = \mathop {\max }\limits_a \,\frac{{{a^T}{\bf{A}}a}}{{{{\left\| a \right\|}^2}}},$ and hence by \eqref{eq: aAa} we have $\lambda_1\leq C.$
Now, since $p/n^2\rightarrow 0$ by assumption, it follows that
\begin{equation}\label{eq: forb_norm_A}
  \frac{{\left\| {\bf{A}} \right\|_F^2}}{{{n^2}}} = \frac{{trace\left( {{{\bf{A}}^2}} \right)}}{{{n^2}}} = \frac{{\sum\limits_{j = 1}^p {\lambda_{j{\bf A}}^2} }}{{{n^2}}} \le \frac{{p\lambda_{1{\bf A}}^2}}{{{n^2}}} \le \frac{{p{C^2}}}{{{n^2}}} \to 0.  
\end{equation}
2) We show that if $p/n^2\rightarrow 0 $, $E(Y^4), E(X_{ij}^4) \leq C$ for $j=1,...,p$ and $C\geq 1$,  and the columns of $\bf{X}$ are independent, then $\frac{\|{\bf A}\|_F^{2}}{n^2}\rightarrow 0$. For $a \in \mathbb{R}^p$ we have by Cauchy–Schwarz,
\begin{align*}
    {a^T}{\bf{A}}a &= E\left[ {Y_i^2{{\left( {\sum\limits_{j = 1}^p {{a_j}{X_{ij}}} } \right)}^2}} \right] \le {\left[ {E\left( {Y_i^4} \right)} \right]^{1/2}}{\left\{ {E\left[ {{{\left( {\sum\limits_{j = 1}^p {{a_j}{X_{ij}}} } \right)}^4}} \right]} \right\}^{1/2}}\\
    &\le {C^{1/2}}{\left\{ {\sum\limits_{{j_1}{j_2}{j_3}{j_4}}^{} {{a_{{j_1}}}{a_{{j_2}}}{a_{{j_3}}}{a_{{j_4}}}E\left( {{X_{i{j_1}}}{X_{i{j_2}}}{X_{i{j_3}}}{X_{i{j_4}}}} \right)} } \right\}^{1/2}}\\
    &\le {C^{1/2}}{\left\{ {\sum\limits_{j = 1}^p {a_j^4E\left( {X_{ij}^4} \right) + \sum\limits_{j \ne j'}^{} {a_j^2a_{j'}^2E\left( {X_{ij}^2} \right)E\left( {X_{ij'}^2} \right)} } } \right\}^{1/2}}\\
    &\le {C^{1/2}}{\left\{ {C\sum\limits_{j = 1}^p {a_j^4 + C\sum\limits_{j \ne j'}^{} {a_j^2a_{j'}^2} } } \right\}^{1/2}}\\
    &\le C{\left\{ {\sum\limits_{j = 1}^p {a_j^4 + \sum\limits_{j \ne j'}^{} {a_j^2a_{j'}^2} } } \right\}^{1/2}} = C{\left\| a \right\|^2}.
\end{align*}
Notice that since that the columns of $\bf{X}$ are independent, the expectation $E\left( {{X_{i{j_1}}}{X_{i{j_2}}}{X_{i{j_3}}}{X_{i{j_4}}}} \right)$ is not zero (up to permutation) when $j_1=j_2$ and $j_3 = j_4$ or when $j_1=j_2=j_3=j_4.$
Also notice we obtained the same result as in \eqref{eq: aAa}, and hence $\frac{\|{\bf A}\|_F^{2}}{n^2}\rightarrow 0$ follows by the same arguments as in the previous example. 
\end{remark}

\begin{remark}\label{c_star_single}
\textbf{\emph{Calculations for Equation \ref{eq:c_star_g}}}:\\ 
Write,
\begin{align*}
\cov \left( {{{\hat \tau }^2},Z_g} \right) &= \cov\left( {\frac{1}{{n\left( {n - 1} \right)}}\sum\limits_{{i_1} \ne {i_2}}^{} {\sum\limits_{j = 1}^p {{W_{{i_1}j}}{W_{{i_2}j}}} } ,\frac{1}{n}\sum\limits_{i = 1}^n {g(X_i)} } \right)\\
&= \frac{1}{{{n^2}\left( {n - 1} \right)}}\sum\limits_{{i_1} \ne {i_2}}^{} {\sum\limits_{j = 1}^p {\sum\limits_{i = 1}^n {E\left( {{W_{{i_1}j}}{W_{{i_2}j}}g(X_i)} \right)} } } \\
&= \frac{2}{{{n^2}\left( {n - 1} \right)}}\sum\limits_{{i_1} \ne {i_2}}^{} {\sum\limits_{j = 1}^p {E\left( {{W_{{i_1}j}}{g_{{i_1}}}} \right)E\left( {{W_{{i_2}j}}} \right)} } \\
&= \frac{2}{{{n^2}\left( {n - 1} \right)}}\sum\limits_{{i_1} \ne {i_2}}^{} {\sum\limits_{j = 1}^p {E\left( {{W_{{i_1}j}}{g_{{i_1}}}} \right){\beta _j}} } \\ 
&= \frac{2}{n}\beta^T\theta_g ,
 \end{align*}
where  \({\theta _g} = E\left[ {Wg\left( X \right)} \right]\).  Also notice that 
\(\var\left( Z_g \right) = \var \left( {\frac{1}{n}\sum\limits_{i = 1}^n {g(X_i)} } \right) = \frac{{\var \left( {g(X)} \right)}}{n}.\)
Thus, by  \eqref{eq:c_star} we get 
$${c_g^*} = \frac{{{\cov} \left( {{{\hat \tau }^2},Z_g} \right)}}{{\var \left( {Z_g} \right)}} = \frac{2\beta^T\theta_g}{{\var \left( {g(X)} \right)}}.$$
\end{remark}

\begin{proof}[Proof of Proposition \ref{prop:singel_asymptotic}]$ $\newline
We wish to prove that  $\sqrt n \left[    T_h  -   T_{\hat h}  \right]\overset{p}{\rightarrow}~0.$
Write,
$$\sqrt n \left( {{T_h} - {T_{\hat h}}} \right) = \sqrt n \left[ {\left( {{T_h} - {T_f}} \right) +  
\left( {{T_f} - {T_{\hat f}}} \right) +
\left( {{T_{\hat f}} - {T_{\hat h}}} \right)}    \right].$$ 
Thus,  we need to show that
 \begin{align}
 & \sqrt n \left( { T_h -  T_f} \right) \overset{p}{\rightarrow} 0, \label{eq:1st_to_show}\\ 
  & \sqrt n \left( { T_{\hat f} -  T_{\hat h}} \right) \overset{p}{\rightarrow}~0, \text{and}\label{eq:2st_to_show}\\
  & \sqrt n \left( { T_{f} -  T_{\hat f}} \right) \overset{p}{\rightarrow} 0, \label{eq:3rd_to_show}
 \end{align}
The proofs of  \eqref{eq:1st_to_show} and \eqref{eq:2st_to_show}  
are essentially the same as the proof of Proposition 5 in \citet{livne2021improved}. This is true since
  $ h(X)\mathbbm{1}_A = f(X)\mathbbm{1}_A,$ where 
        $A$ denotes the event that the selection procedure   $\delta$ perfectly selects the  set $\Theta$, i.e., \(A \equiv \left\{ {{{\bf{S}}_\delta} = \Theta} \right\},\)  and  $\mathbbm{1}_A$ denotes the indicator of   $A.$

We now wish to prove \eqref{eq:3rd_to_show}.
Write,
\[\sqrt n \left[    T_{\hat f}  -  T_f   \right] = \sqrt n \left[ {{{\hat \tau }^2} - {\hat c_f^*}Z_f - \left( {{{\hat \tau }^2} - {c_f^*}Z_f} \right)} \right] =   {\sqrt n }Z_f  {\left( {{{ c_f^*}} - {\hat c_f^*}} \right)} .\]  
By Markov and Cauchy-Schwarz inequalities, it is enough to show that
$$P\left\{ {\left| {\sqrt n Z_f\left( {\hat c_f^* - c_f^*} \right)} \right| > \varepsilon } \right\} \le \frac{{E\left\{ {\left| {\sqrt n Z_f\left( {\hat c_f^* - c_f^*} \right)} \right|} \right\}}}{\varepsilon } \le \frac{{\sqrt {nE\left( {Z_f^2} \right)E\left[ {{{\left( {\hat c_f^* - c_f^*} \right)}^2}} \right]} }}{\varepsilon }\rightarrow 0. $$
Notice that $E(Z_f^2)=\var(Z_f)= \frac{\var[f(X)]}{n}$ and $E[(\hat c_f^*- c_f^*)^2]=\var(\hat c_f^*),$
where by \eqref{eq:list_def} we have
$$
\hat c_f^* = \frac{{\frac{2}{{n\left( {n - 1} \right)}}\sum\limits_{{i_1} \ne {i_2}}^{} {W_{{i_1}}^T{W_{{i_2}}}f\left( {{X_{{i_2}}}} \right)} }}{{\var\left[ {f\left( X \right)} \right]}}\equiv\frac{U}{\var[f(X)]}.
$$
Hence,
it enough to show that
\begin{equation}\label{eq: need to show}
\frac{{\var \left( U \right)}}{{\var [f(X)]}}\rightarrow 0.    
\end{equation}

The variance of $U$ is
\begin{align}\label{eq:Var_U}
\var(U) &= \var\left[ {\frac{2}{{n\left( {n - 1} \right)}}\sum\limits_{{i_1} \ne {i_2}}^{} {\sum\limits_{j = 1}^p {{W_{{i_1}j}}{W_{{i_2}j}}f(X_{i_2})} } } \right] \nonumber\\
 &= \frac{4}{{{n^2}{{\left( {n - 1} \right)}^2}}}\sum\limits_{j,j'}^p {\sum\limits_{{i_1} \ne {i_2},{i_3} \ne {i_4}}^{} {\cov\left[ {{W_{{i_1}j}}{W_{{i_2}j}}f(X_{i_2}),{W_{{i_3}j'}}{W_{{i_4}j'}}f(X_{i_4})} \right].} }     
\end{align}
The covariance in \eqref{eq:Var_U} is different from zero in the following two  cases:
\begin{enumerate}
    \item When $\left\{ {{i_1},{i_2}} \right\}$ is equal to $\left\{ {{i_3},{i_4}} \right\}.$
     \item When one of $\left\{ {{i_1},{i_2}} \right\}$ equals to $\left\{ {{i_3},{i_4}} \right\}$ while the other is different. 
    \end{enumerate}
  The first condition includes two different sub-cases and each of those consists of  $n(n-1)$  quadruples $(i_1,i_2,i_3,i_4)$   that satisfy the condition.
  Similarly, the second  condition above includes four different sub-cases and each of those consists of  $n(n-1)(n-2)$  quadruples  that satisfy the condition. 
  
  We now calculate the covariance for all these six sub-cases.\newline
(1) The covariance  when $ {i_1} = {i_3},{i_2} = {i_4}$ is
 \begin{align*}
 \cov\left[ {{W_{i_1j}}{W_{i_2j}} f(X_{i_2}),{W_{i_1j'}}{{W_{i_2j'}}} f(X_{i_2})} \right] &= E\left( {{W_{i_1j}}{W_{i_1j'}}} \right)E\left[ {{W_{i_2j}}{{ W_{i_2j'}}} f^2(X_{i_2})} \right] \\
 &- E( {{W_{i_1j}}} )E\left[ {{W_{i_2j}} f(X_{i_2})} \right]E\left( {{{ W_{i_1j'}}}} \right)E\left[ {{{ W_{i_2j'}}} f(X_{i_2})} \right]\nonumber\\
 &= E\left( {{W_{i_1j}}{W_{i_1j'}}} \right)E\left[ {{W_{i_2j}}{{ W_{i_2j'}}} f^2(X_{i_2})} \right] - {\beta _j}{\beta _{j'}}{\theta_{fj}}{\theta_{fj'}},
\end{align*}
where recall that $b\equiv Wf(X).$
 Thus, we define
 \begin{equation}\label{eq:delta_1}
    {\delta _1} \equiv \sum\limits_{j,j'}^{} {\left\{ {E\left( {{W_{i_1j}}{W_{i_1j'}}} \right)E\left[ {{W_{i_2j}}{W_{i_2j'}}{f^2(X_{i_2})}} \right] - {\beta _j}{\beta _{j'}}{\theta_{fj}}{\theta_{fj'}}} \right\}}  = E\left( {{b^T}{\bf{A}}b} \right) - {\left( {{\beta ^T}\theta_f } \right)^2}. 
 \end{equation}
    (2) The covariance when  	${i_1} = {i_4},{i_2} = {i_3}$ is 
 \begin{align*}
    \cov\left[ {{W_{i_1j}}{W_{i_2}} f(X_{i_2}),{{ W_{i_2j'}}}{W_{i_1j'}}f(X_{i_1})} \right] &= E\left[ {{W_{i_1j}}{W_{i_1j'}}f(X_{i_1})} \right]E\left[ {{W_{i_2j}}{{ W_{i_2j'}}}f(X_{i_2})} \right]\\
    &- E( W_{i_1j} )E\left[ {{{ W_{i_2j}}} f(X_{i_2})} \right]E\left( {{{ W_{i_2j'}}}} \right)E\left[ {{W_{i_1j'}}f(X_{i_1})} \right]\nonumber\\
 &= {\left\{ {E\left[ {{W_{ij}}{W_{ij'}}f(X_{i})} \right]} \right\}^2} - {\beta _j}{\beta _{j'}}{\theta_{fj}}{\theta_{fj'}}.
 \end{align*}
 Thus, we define
 \begin{equation}\label{eq:delta2}
 {\delta _2} \equiv \sum\limits_{j,j'}^{} {\left\{ {{{\left[ {E\left( {{W_{ij}}{W_{ij'}}f\left( X_{i} \right)} \right)} \right]}^2} - {\beta _j}{\beta _{j'}}{\theta_{fj}}{\theta_{fj'}}} \right\}}  = \left\| {\bf{B}} \right\|_F^2 - {\left( {{\beta ^T}\theta_f } \right)^2},    
 \end{equation}
  where ${\bf{B}} \equiv E\left( {W{W^T}f\left( X \right)} \right).$\newline
(3) Similarly,   when $\,{i_1} = {i_3},{i_2} \ne {i_4}$ we have 
\begin{equation}\label{eq:delta3}
{\delta _3} \equiv \sum\limits_{j,j'}^{} {\left\{ {\left[ {E\left( {{W_{ij}}{W_{ij'}}} \right){\theta_{fj}}{\theta_{fj'}}} \right] - {\beta _j}{\beta _{j'}}{\theta_{fj}}{\theta_{fj'}}} \right\}}  = {\theta_f ^T}{\bf{A}}\theta_f  - {\left( {{\beta ^T}\theta_f } \right)^2.}
\end{equation}
(4)-(5) When $i_1=i_4, i_2 \neq i_3$ or when $i_2=i_3, i_1 \neq i_4$ we have
\begin{equation}\label{eq:delta4_delta_5}
{\delta _4}= \delta_5 \equiv \sum\limits_{j,j'}^{} {\left\{ {{\theta_{fj}}{\beta _{j'}}E\left[ {{W_{ij}}{W_{ij'}}f\left( X_i \right)} \right] - {\beta _j}{\beta _{j'}}{\theta_{fj}}{\theta_{fj'}}} \right\}}  = {\beta ^T}{\bf{B}}\theta_f  - {\left( {{\beta ^T}\theta_f } \right)^2.}
\end{equation}
(6) When $i_2 = i_4, i_1 \neq i_3$,
\begin{equation}\label{eq:delta_6}
{\delta _6} \equiv \sum\limits_{j,j'}^{} {\left\{ {{\beta _j}{\beta _{j'}}E\left[ {{W_{ij}}{W_{ij'}}{f^2(X_i)}} \right] - {\beta _j}{\beta _{j'}}{\theta_{fj}}{\theta_{fj'}}} \right\}}  = {\beta ^T}{\bf{C}}\beta  - {\left( {{\beta ^T}\theta_f } \right)^2},    
\end{equation}
where ${\bf{C}} = E\left( {W{W^T}f^2{{\left( X \right)}}} \right).$
Thus, plugging-in \eqref{eq:delta_1} - \eqref{eq:delta_6} into \eqref{eq:Var_U}  gives
\begin{equation}\label{eq:var_U_deltas}
    \var\left( U \right) = 4 {\left\{ {\frac{1}{{n\left( {n - 1} \right)}}\left( {{\delta _1} + {\delta _2}} \right) + \frac{{\left( {n - 2} \right)}}{{n\left( {n - 1} \right)}}\left( {{\delta _3} + {\delta _4} + {\delta _5} + {\delta _6}} \right)} \right\}}. 
\end{equation}
 Recall that we wish to show that $ \frac{\var(U)}{\var[f(X)]} \rightarrow 0.$
  Thus,   
   it is enough to show that 
  \begin{equation}\label{eq:delta_1-2}
      {\frac{{\left( {{\delta _1} + {\delta _2}} \right)}}{{{n^2 \var[f(X)]}}}}  \to 0,
  \end{equation} 
    and 
  \begin{equation}\label{eq:delta_3-6}
   {\frac{{\left( {{\delta _3} + {\delta _4} + {\delta _5} + {\delta _6}} \right)}}{{{n\var[f(X)]}}}}  \to 0.    
  \end{equation}
 Consider $\delta_1.$ 
For any square matrix ${\bf M},$ we denote ${\lambda _{1{\bf{M}}}}$ to be the largest eigenvalue of ${\bf M}.$ 
 Write,
$$\frac{{{\delta _1}}}{{{n^2}\var \left[ {f\left( X \right)} \right]}} \le \frac{{E\left( {{b^T}{\bf{A}}b} \right)}}{{{n^2}\var[f(X)]}} \le \frac{{{\lambda _{1{\bf{A}}}}}}{n}\frac{{E\left( {{{\left\| b \right\|}^2}} \right)}}{{n\var[f(X)]}} \le \sqrt {\frac{{\left\| {\bf{A}} \right\|_F^2}}{{{n^2}}}} \frac{{E\left( {{{\left\| b \right\|}^2}} \right)}}{{n\var[f(X)]}} \to 0,$$
where the last inequality holds since we assume that $\frac{{\left\| {\bf{A}} \right\|_F^2}}{{{n^2}}} \to 0$ 
and  $\frac{{E\left( {{{\left\| b \right\|}^2}} \right)}}{{n\var[f(X)]}}$ is bounded.
Similarly,  
  $$
  \frac{{{\delta _2}}}{{{n^2}\var[f(X)]}} \le \frac{{\left\| {\bf{B}} \right\|_F^2}}{{{n^2}\var[f(X)]}} \to 0, 
  $$
where
$\frac{{\left\| {\bf{B}} \right\|_F^2}}{{{n^2\var[f(X)]}}} \to 0$ by assumption.
\newline
Consider now $\delta_3.$ Write,
\[\frac{{{\delta _3}}}{{n\var[f(X)]}} \le \frac{{{\theta_f ^T}{\bf{A}}\theta_f }}{{n\var[f(X)]}} \le \frac{{{\lambda _{1{\bf{A}}}}}}{n}\frac{{{{\left\| \theta_f  \right\|}^2}}}{{\var[f(X)]}} \le  
\sqrt {\frac{{\left\| {\bf{A}} \right\|_F^2}}{{{n^2}}}}\frac{{{{\left\| \theta_f  \right\|}^2}}}{{\var[f(X)]}} \to 0,\]
where the above holds since we assume that $\frac{\|\theta_f\|^2}{\var[f(X)]}$ is bounded  and $\frac{{\left\| {\bf{A}} \right\|_F^2}}{{{n^2}}} \to 0.$ 

Recall that $\delta_4=\delta_5.$  
By the Cauchy–Schwarz inequality and the inequality of arithmetic means we have
\begin{equation}\label{eq:bound_delta3}
 \frac{{{\delta _4}}}{{n\var[f(X)]}} \le \frac{{{\beta ^T}{\bf{B}}\theta_f }}{{n\var[f(X)]}} \le  \frac{{\left| {{\lambda _{1{\bf{B}}}}} \right|}}{n}\frac{{\left\| \beta  \right\|\left\| \theta_f  \right\|}}{{\var[f(X)]}} \le \sqrt{\frac{{\left\| {\bf{B}} \right\|_F^2}}{{{n^2}}}} \frac{{\left( {{{\left\| \beta  \right\|}^2} + {{\left\| \theta_f  \right\|}^2}} \right)}}{{2n\var[f(X)]}}\rightarrow 0,
\end{equation}
where the expression  in \eqref{eq:bound_delta3} converges to zero by similar same arguments as shown above.
\newline
Similarly for  $\delta_6,$
$$
\frac{{{\delta _6}}}{{n\var[f(X)]}} \le \frac{{{\beta ^T}{\bf{C}}\beta }}{{n\var[f(X)]}} \le \frac{{{\lambda _{1{\bf{C}}}}}}{n}\frac{{{{\left\| \beta  \right\|}^2}}}{{\var[f(X)]}} \le 
\sqrt{\frac{{\left\| {\bf{C}} \right\|_F^2}}{{{\{n\var[f(X)]\}^2}}}}
\| \beta\|^2 \to 0,
$$
where $\frac{{\left\| {\bf{C}} \right\|_F^2}}{{\{n\var[f(X)]\}^2}} \to 0$  and $\tau^2\equiv\|\beta\|^2= O(1)$ by assumptions.
Hence, \eqref{eq: need to show} follows.
This completes the  proof that
 $\sqrt n \left[    T_{\hat h}  -  T_h   \right]\overset{p}{\rightarrow}~0.$  
\end{proof}

 \begin{remark}\label{remark:E_norm2_b}
 We first show that if $Y^2$ and $X_{ij}^2$ are bounded for all $j=1,...,p$, $i=1,...,n$ and $p/n=O(1)$, then  $\frac{{E\left( {{{\left\| b \right\|}^2}} \right)}}{{n\var[f(X)]}}$ is bounded.  
 Let $C$ be the upper bound of the maximum of $Y^2$ and $X_{ij}^2$, for $j=1,...,p$ and $i=1,...,n.$ 
 Then,
 $$\frac{{E\left( {{{\left\| b \right\|}^2}} \right)}}{{n\var\left[ {f\left( X \right)} \right]}} = \frac{{\sum\limits_{j = 1}^p {E\left( {b_j^2} \right)} }}{{nE\left[ {{f^2}\left( X \right)} \right]}} = \frac{{\sum\limits_{j = 1}^p {E\left[ {X_{1j}^2{Y^2}{f^2}\left( {{X}} \right)} \right]} }}{{nE\left[ {{f^2}\left( X \right)} \right]}} \le \frac{{{C^2}pE\left[ {{f^2}\left( X \right)} \right]}}{{nE\left[ {{f^2}\left( X \right)} \right]}} = \frac{{{C^2}p}}{n}=O(1),$$
where notice that we used the assumption that  $p/n=O(1).$ 

We now show that under the same assumptions as above, together with the assumptions that $\Theta$ is bounded and  $\var[f(X)]\ge c>0$, then 
 $\frac{{\left\| {\bf{B}} \right\|_F^2}}{{{n^2\var[f(X)]}}}   \to 0.$
Recall that $f\left( {{X}} \right) = \sum\limits_{j < j'}^{} {{X_{ij}}{X_{ij'}}}.$
 Notice that when $\Theta$ is bounded,  and when
    the covariates $X_{ij}$,  $i=1,...,n$, $j=1,...,p$  are bounded,
 then so is  $f(X)$. Let $C$ be the upper bound of $|f(X)|$.
Similarly to \eqref{eq: aAa},  for $a \in \mathbb{R}^p$, we have
 \begin{align}\label{eq: aBA}
 {a^T}{\bf{B}}a &= E\left[ {{a^T}W{W^T}af\left( X \right)} \right] \le E\left[ {\left| {f\left( X \right)} \right|{{\left( {{a^T}W} \right)}^2}} \right] \le CE\left[ {{{\left( {{a^T}XY} \right)}^2}} \right]\nonumber\\
 &= CE\left[ {{Y^2}{{\left( {\sum\limits_{j = 1}^p {{a_j}{X_j}} } \right)}^2}} \right] \le {C^2}E\left[ {{{\left( {\sum\limits_{j = 1}^p {{a_j}{X_j}} } \right)}^2}} \right]\nonumber\\
 &= {C^2}\sum\limits_{j \ne j'}^{} {{a_j}{a_{j'}}\underbrace {E\left( {{X_j}{X_{j'}}} \right)}_0}  + {C^2}\sum\limits_{j = 1}^p {a_j^2E\left( {X_j^2} \right)}  = {C^2}{\left\| a \right\|^2}.
\end{align}
It follows that  
 ${\lambda _{1{\bf{B}}}}  \le \frac{{{C^2}{{\left\| a \right\|}^2}}}{{{{\left\| a \right\|}^2}}} = {C^2},$
and by a similar argument as in \eqref{eq: forb_norm_A} we conclude that 
 $$\frac{{\left\| {\bf{B}} \right\|_F^2}}{{{n^2\var[f(X)]}}} \le \frac{{p\lambda _{1{\bf{B}}}^2}}{{{n^2\var[f(X)]}}} \le \frac{pC^4}{{{n^2c}}} \to 0,$$
where recall we assume that $n/p=O(1)$.
A similar argument can be used to show that under the above conditions, 
$\frac{{\left\| {\bf{C}} \right\|_F^2}}{{{\{n\var[f(X)]\}^2}}} \to 0.$
\end{remark}

 \begin{remark}\label{selection_algorithm}
We use the  the following simple selection algorithm $\delta$:

\begin{algorithm}[H]\label{alg1}
\SetAlgoLined
\vspace{0.4 cm}

 \textbf{Input:}
 A dataset \(\left( {{{\bf{X}}_{n \times p}},{{\bf{Y}}_{n \times 1}}} \right)\).
\begin{enumerate}
  \item Calculate $\hat\beta_1^2,...,\hat\beta_p^2$ where  $\hat\beta_j^2$  is given in (\ref{beta_j_hat}) for $j=1,...,p.$   
  
  \item Calculate the  differences
  \({\lambda_{j}} = \hat \beta _{\left( j \right)}^2 - \hat \beta _{\left( {j - 1} \right)}^2\) for $j=2,\ldots,p$ where \(\hat \beta _{\left( 1 \right)}^2 < \hat \beta _{\left( 2 \right)}^2 < ... < \hat \beta _{\left( p \right)}^2\) denotes the  order statistics. 
  \item Select the covariates  \({{\bf{S}}_\delta} = \left\{ {j:\hat \beta _{\left( j \right)}^2 > \hat \beta _{\left( {{j^*}} \right)}^2} \right\}\),  where \({j^*} = \mathop {\arg \max }\limits_j {\lambda_{j}}\).
  \end{enumerate}
\KwResult{Return ${{\bf{S}}_\delta}$. }
  \caption{Covariate selection $\delta$}
\end{algorithm}
The algorithm above finds the largest gap between the ordered estimated squared coefficients and then uses this gap as a threshold to select a set of  coefficients ${{\bf{S}}_\delta} \subset \left\{ {1,...,p} \right\}.$ The algorithm  works well in  scenarios where a relatively large gap truly separates  between  larger coefficients  and the  smaller coefficients of the vector $\beta$.
 \end{remark}

\newpage
\vskip .65cm
\noindent
Ilan Livne (ilan.livne@campus.technion.ac.il)

\noindent
David Azriel (davidazr@technion.ac.il)

\noindent
Yair Goldberg (yairgo@technion.ac.il) 

\vskip 2pt

\noindent
The Faculty of Industrial Engineering and Management, Technion.
\end{document}